\documentclass[12pt]{article}
\usepackage{framed}
\usepackage{amsmath}
\usepackage{amsthm}
\usepackage{amsfonts}
\usepackage[utf8]{inputenc}
\usepackage{color}
\usepackage[T1]{fontenc}
\usepackage[margin=2.5cm]{geometry}
\usepackage{stmaryrd}
\usepackage{graphicx}
\usepackage{dsfont}
\usepackage{mathrsfs}
\usepackage{combelow}
\usepackage[font=small,labelfont=bf]{caption}
\usepackage{tikz}
\usetikzlibrary{patterns}

\theoremstyle{plain}

\newtheorem{theo}{Theorem}[section]
\newtheorem{lemme}[theo]{Lemma}
\newtheorem{prop}[theo]{Proposition}
\newtheorem{coro}[theo]{Corollary}

\renewcommand{\sp}{\hspace{0.2cm}}

\newcommand{\dd}{\mathrm{d}}
\newcommand{\R}{\mathbb{R}}
\newcommand{\C}{\mathbb{C}}

\newcommand{\tend}[2]{\underset{#1\to #2}{\longrightarrow} }
\newcommand{\ds}{\displaystyle}

\newcommand{\eps}{\varepsilon}

\newcommand{\N}{\mathbb{N}}

\newcommand{\vs}{\vspace{1mm}}
\newcommand{\der}[2]{\frac{\dd #1}{\dd #2}}

\newcommand{\tgm}{\tilde{\gamma}}

\newcommand{\diam}{\hspace{0.2mm}\mathrm{diam}\hspace{1mm}}
\newcommand{\curl}{\hspace{0.2mm}\mathrm{curl}\hspace{1mm}}

\renewcommand{\tilde}{\widetilde}

\newcommand{\vertiii}[1]{{\left\vert\kern-0.25ex\left\vert\kern-0.25ex\left\vert #1 
    \right\vert\kern-0.25ex\right\vert\kern-0.25ex\right\vert}} %La norme triple, cf ligne du dessous

\title{2D point vortex dynamics in bounded domains: global existence for almost every initial data}
\author{Martin Donati}
\date{}

\begin{document}

\maketitle

\begin{abstract}
    In this paper, we prove that in bounded planar domains with $C^{2,\alpha}$ boundary, for almost every initial condition in the sense of the Lebesgue measure, the point vortex system has a global solution, meaning that there is no collision  between two point-vortices or with the boundary. This extends the work previously done in \cite{marchioro1984vortex} for the unit disk. The proof requires the construction of a regularized dynamics that approximates the real dynamics and some strong inequalities for the Green's function of the domain. In this paper, we make extensive use of the estimates given in \cite{gustafsson1979motion}. We establish our relevant inequalities first in simply connected domains using conformal maps, then in multiply connected domains.
\end{abstract}

\section{Introduction}

Let us begin by recalling the Euler equations for two dimensional incompressible and inviscid fluids. Let $\Omega$ be an open, bounded and connected subset of $\R^2$. We denote by 
 \begin{equation*}
    u : \begin{cases} \Omega \times \R_+ & \rightarrow \R^2 \\ (x,t) &\mapsto u(x,t),\end{cases}
\end{equation*}  
the velocity of a perfect incompressible fluid filling $\Omega$. Then $u$ must verify  the incompressible Euler equations:
\begin{equation*}
    \begin{cases} \partial_t u(x,t) + u(x,t)\cdot \nabla u(x,t) = - \nabla p(x,t), & \forall(x,t) \in \Omega\times\R^*_+ \vs \\  u(x,0) = u_0(x), & \forall x \in \Omega \\ \nabla \cdot u (x,t) = 0,& \forall(x,t) \in \Omega\times\R_+ \\ u(x,t)\cdot n_\Omega(x) = 0,&\forall (x,t) \in \partial\Omega\times\R_+,\end{cases}
\end{equation*}
where $p$ is the pressure within the fluid, $n_\Omega$ is the exterior normal unit vector to $\partial\Omega$ and $u_0$ is the initial velocity at the time $t=0$.
Introducing  the vorticity $\omega = \curl u = \partial_1 u_2 - \partial_2 u_1$, the first equation of Euler's system gives the following equation for the vorticity:
\begin{equation}\label{eqvorticite}
    \partial_t \omega(x,t) + u(x,t)\cdot\nabla \omega(x,t) = 0.
\end{equation}
We can also express the velocity in terms of the vorticity thanks to the Biot-Savart law. When $\Omega$ is simply connected, the Biot-Savart law reads
\begin{equation}\label{BSG}
    u(x,t) =  \int_\Omega \nabla^{\bot}_x G_\Omega(x,y) \omega(y,t) \dd y,
\end{equation}
where $G_\Omega$ is the Green's function of the domain $\Omega$. It is important to recall that the Green's function of any open and connected subset of $\R^2$ - we call them \emph{domains} in this paper - can be decomposed as
\begin{equation}\label{greendecomp}
    G_\Omega= G_{\R^2} + \gamma_\Omega,
\end{equation}
where $\gamma_\Omega$ is a smooth function in $\Omega\times\Omega$. Indeed, $\gamma_\Omega(x,y)$ is harmonic in both variables.

We define the point vortex system as in \cite{marchioro1993mathematical}. We assume that at the initial time, the vorticity is a sum of Dirac masses $ \omega_0 = \sum_{i=1}^N a_i \delta_{x_i}$ where $N$ is an integer greater than 1, which we fix for the rest of this paper, and the masses $a_i$ are real numbers, also fixed. Since the vorticity equation \eqref{eqvorticite} is a transport equation in $\omega$, we expect the vorticity to remain a sum of Dirac masses with the same intensity as at the initial time. So we choose to write $ \omega(t) = \sum_{i=1}^N a_i \delta_{x_i(t)}$. We then define the point vortex system in a simply connected domain as the solution of the system of equations:
\begin{equation}
\label{ptvortexdynamic}    \forall 1 \le i \le N, \sp \der{x_i(t)}{t} = \sum_{\substack{j= 1\\ j\neq i}}^N \nabla_x^\bot G_\Omega(x_i(t),x_j(t))a_j + \nabla_x^\bot\gamma_\Omega(x_i(t),x_i(t)) a_i.
\end{equation}
This is obtained by introducing the expression of the vorticity and the decomposition \eqref{greendecomp} into the Biot-Savart law \eqref{BSG} and by removing the singular term that appears in the limit of $ \nabla_x G_{\R^2}(x,x_i(t)) = \frac{(x-x_i(t))^\bot}{2\pi|x-x_i(t)|^2}$ when $x$ goes to $x_i(t)$. This term represents  high speed rotation around $x_i(t)$, so it shouldn't affect the motion of $x_i(t)$ itself. 

All those choices have been mathematically justified in \cite{marchioro1993VorticiesAndLocalization}, where it has been proved that highly concentrated smooth solution of the Euler equations converges in the sense of measures to the solution of the point vortex system, as the initial data converges towards the initial sum of Dirac masses. It has also been proved in \cite{ConvergencePVtoEulerGoodman} that the point vortex system is a good approximation of the Euler equation from the point of view of numerics, taking as initial data a grid of vortices approaching a smooth initial vorticity. 

The question that naturally arises now is whether the system of equations defined in \eqref{ptvortexdynamic} has a global solution for every initial condition $(x_i(0))_i$. The answer to that question unfortunately is negative in general, since in $\R^2$ one can build an initial datum such that point vortices collapse in finite time. See \cite{marchioro1993mathematical}, \cite{iftimielargetime}, or \cite{IftimieMarchioroSelfSimilar} for explicit examples. By construction, the point vortex dynamics isn't defined anymore as soon as a collapse occurs, since equation \eqref{ptvortexdynamic} becomes singular when two points collide. But what we can expect is that these occurrences of collapse are exceptional, meaning that the initial configurations leading to collapse are negligible in the sense of the Lebesgue measure. This result has been proved in \cite{marchioro1984vortex} in the unit disk $D(0,1)$. In the case of $\R^2$, it has been proved with the additional assumption that every possible sum of the masses never vanishes, meaning that $\sum_{i\in P} a_i \neq 0$ for every $P \subset \{1,\ldots,n\}$. Proofs of these results can be found in \cite{marchioro1984vortex} and \cite{marchioro1993mathematical}. Very recently, \cite{godardcadillac2021vortex} proved that the assumption that $\sum_{i=1}^N a_i \neq 0$ in $\R^2$ can be removed.

Let us give a precise statement of the result of \cite{marchioro1984vortex} in the case of the disk.
\begin{theo}\label{theomarchioro}
If $\Omega = D(0,1)$ the open unit disk, then the point vortex dynamics \eqref{ptvortexdynamic} for any fixed number of points $N \ge 1$ and masses $(a_i)_i \in \R^N$ is globally well defined except maybe for a set of initial conditions in $\Omega^N$ which has vanishing Lebesgue measure.
\end{theo}

The purpose of this article is to prove a generalization of this theorem to more general bounded domains. Let $\Omega$ be an open bounded and connected subset of $\R^2$ with a $C^{2,\alpha}$ boundary for some $\alpha >0$. In the case where $\Omega$ is multiply connected, the point vortex dynamics is changed since the Biot Savart law is different. We refer to \cite[Chapter 15]{flucher_variational_1999} for the point-vortex system in multiply-connected domains. We will give all the details in section \ref{sectionMultiConnected} but the result is that for a domain $\Omega$ that has $m$ holes, the point vortex dynamics is given by 
\begin{equation}\label{ptVortexDynamicCOmplete}
 \der{x_i(t)}{t} = \sum_{\substack{j= 1\\ j\neq i}}^N \nabla_x^\bot G_\Omega(x_i(t),x_j(t))a_j + \nabla_x^\bot\gamma_\Omega(x_i(t),x_i(t)) a_i + \sum_{j=1}^m c_j(t) \nabla^\perp w_j(x_i)
\end{equation}
for all $1 \le i \le N$. Above 
\begin{equation*}
c_j(t)=\xi_j+\sum_{k=1}^N a_k w_j(x_k(t)),
\end{equation*}
$\xi_j$ is the circulation of the velocity $u$ on the boundary of the $j$-th hole of $\Omega$,  and $w_j$ are the harmonic measures of the domain $\Omega$.
Let us observe that by the Kelvin theorem, the circulations $\xi_j$ are constant in time. They are therefore prescribed at the initial time.

Let $\lambda$ be the Lebesgue measure on $\R^{2N}$. We define on the set $\overline{\Omega}^N$:
\begin{equation*}
    d(X) = \min \left( \min_{i\neq j} |x_i - x_j|, \min_i \mathrm{d}(x_i,\partial\Omega) \right)\quad\forall X = (x_1,\ldots,x_N).
\end{equation*}

We define $\Gamma = \{ X = (x_1,\ldots,x_N), d(X) > 0\}. $ This is the set of all configurations for which relation \eqref{ptVortexDynamicCOmplete} makes sense. We note by $S_tX$ the evolved configuration by the dynamics \eqref{ptVortexDynamicCOmplete} from the starting configuration $X\in\Gamma$, after a time $t$. We know there exists a time $\tau(X) = \sup \{t \ge 0, S_tX \in \Gamma\} > 0$ until which the dynamics is well defined. In this paper we will prove the following theorem.

\begin{theo}\label{mainresult}
Let $\Omega$ be an open, bounded and connected subset of $\R^2$ with a $C^{2,\alpha}$ boundary for some $\alpha >0$ with $m\in \N$ holes. We fix the number of point vortices $N\ge 1$, the masses $(a_i)_{1\le i \le N} \in \R^N$, and the circulations $(\xi_j)_{1\le j \le m} \in \R^m$.  With the previous notations we have that
\begin{equation*}
    \lambda (\{ X\in \Omega^N, \tau(X) < \infty\} ) = 0,
\end{equation*}
meaning that for almost every starting position $X$, the point vortex dynamics in $\Omega$ is well defined for every time.
\end{theo}

We observe that Theorem \ref{mainresult} has a simple proof in the case of a single point vortex in a simply connected domains, that is in the case $N=1$ and $m=0$. We introduce the Robin function $\tgm_\Omega(x) = \gamma_\Omega(x,x)$. Since $\gamma_\Omega(x,y) = \gamma_\Omega(y,x)$, we have that $\nabla\tgm_\Omega(x) = \nabla_x\gamma_\Omega(x,x) + \nabla_y \gamma_\Omega(x,x) = 2 \nabla_x\gamma_\Omega(x,x)$. In this case the dynamics of a single point vortex becomes
\begin{equation*}
    \der{x(t)}{t} = \frac{1}{2}\nabla^\bot\tgm_\Omega(x(t))a.
\end{equation*}
Therefore, a single point vortex evolves on the level set of the Robin function. This map has been studied in \cite{gustafsson1979motion}, from which we know in particular that $\tgm_\Omega(x) \tend{x}{\partial\Omega} + \infty$.  Therefore a single point vortex can't hit the boundary so the dynamics is well defined for every time.

The proof of Theorem \ref{mainresult} that we will give in section \ref{sectionMainResult} borrows arguments from \cite{marchioro1984vortex}, but we have to deal with two major difficulties. The first one is the construction of a convenient regularized dynamics, and the second one is to prove some analytic inequalities on the Green's function and the Robin function of the domain $\Omega$. In Section \ref{sectionConformal} we will obtain the required inequalities for those maps in the case of simply connected domains, and in Section \ref{sectionMultiConnected}, we will extend those results to the case of multiply connected domains. Section \ref{sectionPreuve} is devoted to the construction of the regularized dynamics and the completion of the proof of Theorem \ref{mainresult}.

\section{Simply connected and exterior domains}\label{sectionConformal}

\small
List of notations:
\begin{itemize}
    \item $N\in \N$ denotes the number of point vortices;
    \item $\lambda$ is the Lebesgue measure on $\R^{2N}$;
    \item $(x_1,x_2)^\perp = (-x_2,x_1)$;
    \item $\Omega$ is a $C^{2,\alpha}$ bounded domain of $\R^2$ with $m \in \N$ holes, and its boundaries are $\Gamma_j$, $0 \le j \le m$, with $\Gamma_0$ the exterior boundary;
    \item $\mathcal{U}$ denotes a general bounded domain with $C^{2,\alpha}$ boundary;
    \item $U$ denotes a general simply connected bounded domain with $C^{2,\alpha}$ boundary;
    \item $\Pi$ denotes a general exterior domain with $C^{2,\alpha}$ boundary;
    \item $D(x_0,r)$ is the disk of center $x_0$ and of radius $r$ and $D = D(0,1)$;
    \item $\Pi_D = \big(\overline{D(0,1)}\big)^c$;
    \item $T$ denotes a biholomorphic map, usually from $U$ to $D$ or from $\Pi$ to $\Pi_D$;
    \item $n_\mathcal{U}$ is the exterior normal unit vector to $\partial\mathcal{U}$, extended to a neighborhood of $\partial\mathcal{U}$ by relation \eqref{eqTranspVecteurNormal} when possible;
    \item $G_\mathcal{U}$ is the Green's function of the domain $\mathcal{U}$, and $G = G_\Omega$ in section \ref{sectionPreuve};
    \item $\gamma_\mathcal{U}$ is the regular part of $G_\mathcal{U}$, see relation \eqref{greendecomp}, and $\gamma = \gamma_\Omega$ in section \ref{sectionPreuve};
    \item $\tgm_\mathcal{U}(x) = \gamma_\mathcal{U}(x,x)$ is the Robin function of the domain $\mathcal{U}$, and $\tgm = \tgm_\Omega$ in section \ref{sectionPreuve} ;
    \item $C, C_1, C_2, \ldots$, are strictly positive constants that may vary from one line to another, when their value is not important to the result;
    \item $a \cdot b$ is the scalar product of vectors in $\R^2$;
    \item $\nabla f$ and $\nabla \cdot g$ are respectively the gradient of $f$ and the divergence of $g$;
    \item $V_j$ are neighborhoods of $\Gamma_j$, and $K$ is a compact set as in the decomposition \eqref{decompOmegaVj};
    \item $S_tX$ is the solution of the point vortex dynamics starting from $X$ after a time $t$, and $S_t^\eps X$ the regularized dynamics constructed in Section \ref{sectionRegDynamic};
\end{itemize}
\normalsize 

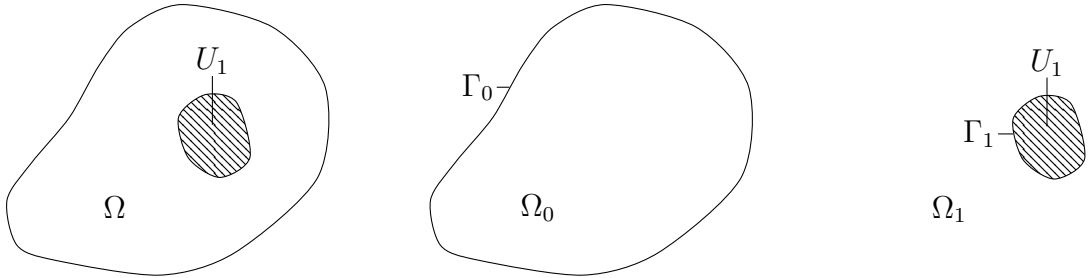
\begin{figure}\centering
    \begin{tikzpicture}
    \draw[fill=white] plot [smooth cycle]coordinates  {(1,1) (2.4,0.8) (3.4, 1.1) (4.5, 2.1) (4.6, 3.3) (3.9,4.1) (2.7,4.4) (2.0, 4.1) (1.6,3.6) (1.2, 2.9) (0.7, 2.3) (0.37,1.8) (0.5,1.2)} node at (1.8,1.7) {$\Omega$};
    \draw[pattern=north west lines]   plot [smooth cycle]coordinates  {(2.65,2.90) (2.78,2.36) (3.2,2.1) (3.6,2.37) (3.42,3.10) (3,3.2)} ;
    \draw node at (3.1,3.65) {$U_1$};
    \draw (3.1,3.45) -- (3.1,2.8);
    \end{tikzpicture}
    \hspace{1cm}
    \begin{tikzpicture}
    \draw plot [smooth cycle]coordinates  {(1,1) (2.4,0.8) (3.4, 1.1) (4.5, 2.1) (4.6, 3.3) (3.9,4.1) (2.7,4.4) (2.0, 4.1) (1.6,3.6) (1.2, 2.9) (0.7, 2.3) (0.37,1.8) (0.5,1.2)} node at (1.8,1.7) {$\Omega_0$}  node at (1,3.3) {$\Gamma_0$};
    \draw(1.25,3.3) --(1.42,3.3);
    \end{tikzpicture}
    \hspace{1cm}
    \begin{tikzpicture}
    \fill[fill = white] plot [rectangle]coordinates  {(0.5,0.8) (4.6,0.8) (4.6,4.4) (0.5,4.4)};
    \draw[pattern=north west lines]  plot [smooth cycle]coordinates  {(2.65,2.90) (2.78,2.36) (3.2,2.1) (3.6,2.37) (3.42,3.10) (3,3.2)} node at (1.8,1.7) {$\Omega_1$} node at (2.2,2.7) {$\Gamma_1$} ;
    \draw(2.45,2.7) --(2.65,2.7);
    \draw node at (3.1,3.65) {$U_1$};
    \draw (3.1,3.45) -- (3.1,2.8);
    \end{tikzpicture}
    \caption{An example for $m=1$. The domains $\Omega$, $\Omega_0$ and $\Omega_1$ and their boundaries.}
    \label{fig1}
\end{figure}

For the rest of this paper, $\Omega$ denotes a bounded domain whose boundary is $C^{2,\alpha}$ for some $0<\alpha<1$. It is either simply connected or it has $m\in \N$ holes that are the simply connected bounded domains $U_1,\ldots,U_m$ and their boundaries are $\Gamma_1,\ldots,\Gamma_m$. We denote by $\Gamma_0$ the exterior boundary of $\Omega$, meaning that $\Omega$ lies within the interior in the Jordan sense of $\Gamma_0$, and by $\Omega_0$ the simply connected bounded domain whose boundary is $\Gamma_0$, namely the domain $\Omega$ "without holes". Finally, we call \emph{exterior domain} a domain whose complement is bounded and simply connected. We denote for $1\le j\le m$,  $\Omega_j = (\overline{U_j})^c$ the exterior domains of the $m$ holes. The domain $\Omega$ is pictured in Figure \ref{fig1}.

Since $\Omega = \bigcap_{j=0}^m \Omega_j$, with $\Omega_0$ simply connected, and where $\Omega_j$, $j\ge 1$ are exterior domains, our strategy in this paper is to establish inequalities on $\Omega$ by establishing them for any bounded and simply connected domain $U$, and for any exterior domain $\Pi$.

We also denote by $D = D(0,1)$ the unit disk, and by $\Pi_D$ the exterior domain of $D$, namely $\Pi_D = \{ x \in \C, |x| >1 \}$. 

\subsection{Holomorphic maps}
Holomorphic maps are the subject of the Chapter 2 of \cite{ahlfors1966complexanalysis}. A map $T : \C \rightarrow \C$ is a \emph{biholomorphic map} if $T$ and $T^{-1}$ are holomorphic maps. Such maps satisfy that their derivative never vanishes. Let us write $T = T_1 + iT_2$, and we identify $\R^2$ and $\C$, meaning that we also denote $ T = \begin{pmatrix} T_1 \\ T_2 \end{pmatrix}$. Then a holomorphic map satisfies the Cauchy-Riemann equations
\begin{equation*}
    \begin{cases}\partial_1 T_1 = \partial_2 T_2 \\ \partial_1 T_2 = - \partial_2 T_1.\end{cases}
\end{equation*} 
We have that
\begin{equation*}
    T' = \partial_1T = -i\partial_2 T
\end{equation*}
and therefore
\begin{equation*}
    T'' = \partial_1^2 T = - \partial_2^2T = -i\partial_1\partial_2T.
\end{equation*}
Finally, the Jacobian matrix of $T$ is $JT = \begin{pmatrix}\partial_1 T_1 & \partial_2 T_1 \\ \partial_1 T_2 & \partial_2T_2 \end{pmatrix}$, so $\det JT = |T'|^2$. In the following, we will freely use these properties. In particular, we can always substitute the second partial derivative $\partial_2$ with $i\partial_1$ according to these formulas. For any map $f \in C^1(\C,\R)$, we have that
\begin{equation}\label{eqfRondT}
    \nabla (f\circ T)(x) = \begin{pmatrix} \partial_1 T_1(x) \partial_1 f(T(x)) + \partial_1 T_2(x) \partial_2 f(T(x)) \\ - \partial_1 T_2(x) \partial_1 f(T(x)) + \partial_1 T_1(x) \partial_2 f(T(x))    \end{pmatrix}.
\end{equation}

We conclude this paragraph with a technical lemma.
\begin{lemme}\label{lemmeIneqIntermediaires}
For any bounded domain $\mathcal{U}$ whose boundary is $C^{2,\alpha}$, and for any $\kappa <1$, we have that
\begin{equation*}
    \iint_{\mathcal{U}\times\mathcal{U}} \frac{1}{|x-y|^{1+\kappa}} \dd x \dd y  < \infty,
\end{equation*}
and
\begin{equation*}
    \iint_{\mathcal{U}\times\mathcal{U}} \frac{1}{d(x,\partial\mathcal{U})^{\kappa}} \frac{1}{|x-y|} \dd x \dd y < \infty.
\end{equation*}
Moreover there exists a constant $C$ depending only on $\mathcal{U}$ such that for sufficiently small $\eps > 0$, 
\begin{equation}\label{star}
    \int_{\{x \in \mathcal{U}, d(x,\partial\mathcal{U}) \ge \eps \}} \frac{1}{d(x,\partial\mathcal{U})^{1+\kappa}}\dd x \le C \eps^{-\kappa}.
\end{equation}
\end{lemme}
\begin{proof}
Let $R=\diam(\mathcal{U})$ so that $\mathcal{U} \subset B(x,R)$ for every $x \in \mathcal{U}$. Then
\begin{equation*}
    \iint_{\mathcal{U}\times\mathcal{U}} \frac{1}{|x-y|^{1+\kappa}} \dd x \dd y  \le \int_\mathcal{U}\int_{B(x,R)} \frac{1}{|x-y|^{1+\kappa}}\dd x \dd y = |\mathcal{U}| \int_0^R\int_0^{2\pi} \frac{1}{r^{1+\kappa}}r\dd \theta\dd r < \infty.
\end{equation*}
With the same argument,
\begin{equation*}
    \iint_{\mathcal{U}\times\mathcal{U}} \frac{1}{d(x,\partial\mathcal{U})^{\kappa}} \frac{1}{|x-y|} \dd x \dd y  \le 2\pi R \int_\mathcal{U}\frac{1}{d(x,\partial\mathcal{U})^{\kappa}} \dd x.
\end{equation*}
To prove that the integral $ \int_\mathcal{U}\frac{1}{d(x,\partial\mathcal{U})^{\kappa}} \dd x$ is finite, and to prove \eqref{star},
we make a finite number of local changes of coordinates and we use that $\mathcal{U}$ is bounded to write
\begin{equation*}
    \int_\mathcal{U}\frac{1}{d(x,\partial\mathcal{U})^{\kappa}} \dd x \le C \int_0^R \frac{1}{s^\kappa}\dd s < \infty
\end{equation*}
and
\begin{equation*}
    \int_{\{x \in \mathcal{U}, d(x,\partial\mathcal{U}) \ge \eps \}} \frac{1}{d(x,\partial\mathcal{U})^{1+\kappa}}\dd x \le \int_\eps^R \frac{1}{s^{1+\kappa}}\dd s \le C \eps^{-\kappa}.
\end{equation*}
\end{proof}

\subsection{The Riemann Mapping Theorem}

We refer now to Chapter 6 of \cite{ahlfors1966complexanalysis}.

\begin{theo}[Riemann Mapping Theorem]
For any non empty, open and simply connected subset $U$ of $\C$, that isn't the whole plane, there exists a biholomorphism from $U$ to the unit disk $D$.
\end{theo}

The consequence of this theorem is that any suitable domain is linked to the disk by a map that has very interesting properties related to the Green's function of both domains. However the Riemann Mapping Theorem only states the theoretical existence of such map, and only a few explicit examples are known. In particular, we have no control over the derivatives of the biholomorphism in general. We combine Theorems 3.5 and 3.6 from \cite{pommerenke1992boundary} to obtain the following corollary of the Kellogg-Warschawski Theorem.

\begin{theo}\label{theoKW}
Let $T$ be a biholomorphism mapping on a bounded, open, and simply connected set $U$ whose boundary  $\partial U$ is a $C^{2,\alpha}$ Jordan curve, with $0<\alpha<1$. Then $T$, $T'$ and $T''$ are continuous up to $\overline{U}$, and $T^{-1}$ and $(T^{-1})'$ are continuous up to $\overline{T(U)}$. Thus there exist constants $m$ and $M$ satisfying for every $x\in\overline{U}$ that $0<m\le|T'(x)|\le M$ and $|T''(x)| \le M$.
\end{theo}

Let us stress the fact that since the automorphisms of the disk are known explicitly and belong to $C^\infty(\overline{D})$, if there exists one biholomorphism $T : U \rightarrow D$ that is smooth up to the boundary, then every biholomorphism $ T: U \rightarrow D$ is smooth up to the boundary. Since in this paper we will always consider smooth domains $U$, every biholomorphism $T : U \mapsto D$ will satisfy the conclusions of Theorem \ref{theoKW}.

Please note that the $C^{2,\alpha}$ condition is not optimal. For instance, it is known that if $\partial U$ has a parametrization with a Dini-continuous curvature, then the conclusion of the theorem is still true. Also, assuming that $\partial U \in C^{n,\alpha}$ implies more generally that $T^{(n)}$ is continuous up to $\partial U$. Despite these remarks, we will stick to the condition $C^{2,\alpha}$ in the context of this article.

In conclusion, the Riemann Mapping Theorem states the existence of the map $T : U \rightarrow D$ and Theorem \ref{theoKW}  states that $\forall x\in\overline{U}, 0<m\le|T'(x)|\le M$ and $|T''(x)| \le M$.

We have a very similar result, this time concerning \emph{exterior domains}, which we define as the complement of the closure of a bounded and simply connected set in $\C$. We have the following theorem.

\begin{theo}\label{theoComplement}
Let $\Pi$ be an exterior domain, with $C^{2,\alpha}$ boundary. Let $T$ be a biholomorphic map from $\Pi$ to $\Pi_D = \{x\in \C, |x| > 1 \}$. Such a map exists and satisfies that $T$, $T'$ and $T''$ are continuous up to $\overline{\Pi}$, and that $T^{-1}$ and $(T^{-1})'$ are continuous up to $\overline{T(\Pi)}$. Moreover  there exist constants $m$ and $M$ such that $\forall x\in\Pi$, $0<m\le|T'(x)|\le M$ and $|T''(x)| \le M$.
\end{theo}
The proof of this result can be found in \cite{IftimieLopesLopesObstacle2003}. It follows from the bounded domain case using the holomorphic map $T : \Pi_D \rightarrow D$,  $T(z) = \frac{1}{z}$.

\subsection{Green's Function}
We start by recalling that for every $(x,y) \in D\times D$, $x \neq y$,
\begin{equation}\label{GreenD}
    G_D(x,y) = \frac{1}{2\pi} \ln \frac{|x-y|}{|x-y^*||y|},
\end{equation}
where $y^* = \frac{y}{|y|^2}$ is the inverse of $y$ relative to the unit circle. Using the decomposition \eqref{greendecomp} we have that for every $(x,y) \in D\times D$, $x \neq y$,
\begin{equation}\label{eqGammaD}
    \gamma_D(x,y) = G_D(x,y) - G_{\R^2}(x,y) = -\frac{1}{2\pi} \ln (|x-y^*||y|)
\end{equation}
and by continuity the relation $\gamma_D(x,y) = -\frac{1}{2\pi} \ln (|x-y^*||y|)$ holds true also for $y=x$. For every $x \in D$ we thus have that
\begin{equation}\label{eqTgmD}
    \tgm_D(x) = \gamma_D(x,x) = -\frac{1}{2\pi} \ln ||x|^2 - 1|.
\end{equation}
Notice that $\tgm_D$ is a radial function.

Let $U$ be a bounded and simply connected domain. If $T : U \rightarrow D$ is a biholomorphic map, then we have the following property.
\begin{prop}\label{transpgreen}
For every $(x,y) \in U\times U$, $x\neq y$, \begin{equation*}
    G_U(x,y) = G_D(T(x),T(y)) = \frac{1}{2\pi} \ln\frac{|T(x)-T(y)|}{|T(x)-T(y)^*||T(y)|}.
\end{equation*}
\end{prop}
Using this in the decomposition \eqref{greendecomp} we obtain
\begin{equation}\label{transpgamma}
    \gamma_U(x,y) + G_{\R^2}(x,y) = \gamma_D(T(x),T(y)) + G_{\R^2}(T(x),T(y))
\end{equation}
and thus
\begin{align*}
     \forall x\in U, \sp \sp \tgm_U(x) & = \lim_{y\rightarrow x} \gamma_U(x,y) \\
    & = \lim_{y\rightarrow x} \left(\gamma_D(T(x),T(y)) + \frac{1}{2\pi}\ln \frac{|T(x) - T(y)|}{|x-y|}\right).
\end{align*}
Therefore
\begin{equation}\label{transptgm}
        \forall x\in U, \sp \sp \tgm_U(x) = \tgm_D(T(x)) + \frac{1}{2\pi} \ln|T'(x)|.
\end{equation}

A quite remarkable fact is that for every $(x,y) \in \Pi_D, x \neq y,$
\begin{equation}\label{GreenPiD}
    G_{\Pi_D}(x,y) = \frac{1}{2\pi} \ln \frac{|x-y|}{|x-y^*||y|},
\end{equation}
which is the same expression as for $G_D(x,y)$. Thus the relations above also hold true for any exterior domain $\Pi$ and any biholomorphism $T : \Pi \rightarrow \Pi_D$, which exists according to Theorem \ref{theoComplement}. For example, for every $(x,y) \in \Pi\times\Pi$, $x\neq y$, \begin{equation*}
    G_\Pi(x,y) = G_{\Pi_D}(T(x),T(y)) = \frac{1}{2\pi} \ln\frac{|T(x)-T(y)|}{|T(x)-T(y)^*||T(y)|}.
\end{equation*}

Let us recall here a classical theorem, see for example \cite[Theorem 4.17]{aubin1982Nonlinear}.
\begin{theo}\label{lemmeGreenReguliere}
Let $\mathcal{U}$ be a bounded domain with $C^{2,\alpha}$ boundary. Then $G_\mathcal{U} \in C^2\big(\overline{\mathcal{U}}\times\overline{\mathcal{U}} \setminus \{ (x,x), x \in \overline{\mathcal{U}}\}\big)$.
\end{theo}
In other words, except where $x=y$, the Green's function is smooth up to the boundary.

In their proof of Theorem \ref{theomarchioro}, the authors of \cite{marchioro1984vortex} show that for any $\kappa <1$,
\begin{equation}\label{eqMarchioro1}
    \iint_{D\times D} \frac{1}{d(x,\partial D)^\kappa}|\nabla_x G_D(x,y) \cdot \nabla^\bot \tgm_D(x)| < \infty
\end{equation}
and
\begin{equation}\label{eqMarchioro2}
    \iint_{D\times D} \frac{1}{|x-y|^\kappa}|\nabla_x G_D(x,y) \cdot \nabla^\bot \tgm_D(x)| < \infty.
\end{equation}

In this paper we will extend these inequalities to the more general bounded domain $\Omega$.

\subsection{Intermediate lemmas}\label{sectionIntermediateLemmas}
We recall that $U$ is a bounded simply connected domain, and $T$ is a biholomorphism from $U$ to $D$.

The following lemma is the first step to extend inequalities \eqref{eqMarchioro1} and \eqref{eqMarchioro2} to the domain $U$.
\begin{lemme}\label{lemmeTransp} For every $x$ and $y$ in $U$:
\begin{equation}\label{transpscalggamma}
    \nabla_x G_U(x,y) \cdot \nabla^\bot \tgm_U(x) = |T'(x)|^2 \nabla_x G_D(T(x),T(y)) \cdot \nabla^\bot\tgm_D(T(x)) + \frac{\nabla_xG_D(T(x),T(y))\cdot \psi(x)}{|T'(x)|^2},
\end{equation}
with $\psi: U \rightarrow \R^2$  an explicit bounded function on $U$.
\end{lemme}

\begin{proof}
From Proposition \ref{transpgreen} and relation \eqref{eqfRondT}, we have that
\begin{equation*}
    \nabla_x G_U(x,y) = \begin{pmatrix} \partial_{x_1} G_D(T(x),T(y)) \partial_1 T_1(x) + \partial_{x_2} G_D(T(x),T(y)) \partial_1 T_2(x) \\ -\partial_{x_1} G_D(T(x),T(y)) \partial_1 T_2(x) + \partial_{x_2} G_D(T(x),T(y)) \partial_1 T_1(x)\end{pmatrix}.
\end{equation*}

Similarly from relation \eqref{transptgm} and relation \eqref{eqfRondT}, since $T'$ is also a holomorphic map, we have that

\begin{equation}\label{eqNablaTgmDev}
    \nabla \tgm_U(x) = \begin{pmatrix}\ds \partial_1 \tgm_D \partial_1 T_1 + \partial_2 \tgm_D  \partial_1 T_2  + \frac{\partial_1 T_1 \partial_1^2 T_1 + \partial_1 T_2 \partial_1^2 T_2}{2\pi|T'|^2} \vspace{2mm}\\ \ds -\partial_1 \tgm_D  \partial_1 T_2 + \partial_2 \tgm_D  \partial_1 T_1  + \frac{-\partial_1 T_1 \partial_1^2 T_2 + \partial_1 T_2 \partial_1^2 T_1}{2\pi|T'|^2}\end{pmatrix}.
\end{equation}
Therefore,
\begin{align*}
    \nabla_x G_U \cdot \nabla^\bot \tgm_U & = - \left( \partial_{x_1} G_D\partial_1T_1 + \partial_{x_2} G_D \partial_1T_2 \right) \\
    & \sp \sp \sp \sp \times \left( -\partial_1 \tgm_D  \partial_1 T_2 + \partial_2 \tgm_D  \partial_1 T_1  + \frac{-\partial_1 T_1 \partial_1^2 T_2 + \partial_1 T_2 \partial_1^2 T_1}{2\pi|T'|^2} \right) \\ 
    & \sp \sp + \left( -\partial_{x_1} G_D\partial_1T_2 + \partial_{x_2} G_D \partial_1T_1 \right) \\
    & \sp \sp \sp \sp \times  \left(  \partial_1 \tgm_D \partial_1 T_1 + \partial_2 \tgm_D  \partial_1 T_2  + \frac{\partial_1 T_1 \partial_1^2 T_1 + \partial_1 T_2 \partial_1^2 T_2}{2\pi|T'|^2}\right). 
\end{align*}
We notice that the terms with the factor $\partial_{x_1} G_D \partial_1 \tgm_D$ cancel each others, as well as the terms with $\partial_{x_2} G_D \partial_2 \tgm_D$. We thus have that
\begin{align*}
    \nabla_x G_U \cdot \nabla^\bot \tgm_U & = - \frac{1}{2\pi}\left( \partial_{x_1} G_D\partial_1T_1 + \partial_{x_2} G_D \partial_1T_2 \right) \left( \frac{-\partial_1 T_1 \partial_1^2 T_2 + \partial_1 T_2 \partial_1^2 T_1}{|T'|^2} \right) \\ 
    & \sp \sp + \frac{1}{2\pi}\left( - \partial_{x_1} G_D\partial_1T_2 + \partial_{x_2} G_D \partial_1T_1 \right) \left( \frac{\partial_1 T_1 \partial_1^2 T_1 + \partial_1 T_2 \partial_1^2 T_2}{|T'|^2}\right) \\
    & \sp \sp - \partial_{x_1} G_D\partial_2 \tgm_D  (\partial_1 T_1)^2 + \partial_{x_2} G_D \partial_1 \tgm_D  (\partial_1 T_2)^2 \\ 
    & \sp \sp - \partial_{x_1} G_D\partial_2 \tgm_D  (\partial_1 T_2)^2 + \partial_{x_2} G_D \partial_1 \tgm_D  (\partial_1 T_1)^2.
\end{align*}

The last two rows can be simplified, showing that they are equal to $|T'|^2 \nabla_x G_D \cdot \nabla^\bot\tgm_D$. For the first two rows, we factor out by $\partial_{x_i} G_D$ and then by $\partial_1^2 T_i$, so that

\begin{align*}
    \nabla_x G_U \cdot \nabla^\bot \tgm_U & = \frac{\partial_{x_1} G_D}{2\pi|T'|^2} [ -2\partial_1T_1\partial_1T_2 \partial_1^2 T_1 + ((\partial_1T_1)^2 - (\partial_1T_2)^2 )\partial_1^2 T_2] \\
    & \sp \sp +  \frac{\partial_{x_2} G_D}{2\pi|T'|^2} [ (-(\partial_1T_2)^2 + (\partial_1T_1)^2) \partial_1^2 T_1 + 2\partial_1T_1\partial_1T_2\partial_1^2 T_2]\\
    & \sp \sp + |T'|^2 \nabla_x G_D \cdot \nabla^\bot\tgm_D.
\end{align*}
This proves equality \eqref{transpscalggamma} with the following explicit function $\psi$:
\begin{equation}\label{DefPsi}
    \psi(x) = \frac{1}{2\pi} \begin{pmatrix}-2\partial_1T_1\partial_1T_2 \partial_1^2 T_1 + ((\partial_1T_1)^2 - (\partial_1T_2)^2 )\partial_1^2 T_2 \\ (-(\partial_1T_2)^2 + (\partial_1T_1)^2) \partial_1^2 T_1 + 2\partial_1T_1\partial_1T_2\partial_1^2 T_2 \end{pmatrix}.
\end{equation}
The function $\psi$ is bounded since, thanks to Theorem \ref{theoKW}, both $T'$ and $T''$ are bounded.
\end{proof}
Notice that  the proof of this lemma only uses Proposition \ref{transpgreen} and relation \eqref{transptgm}, so the lemma also holds in the domain $\Pi$. More precisely, if $T$ denotes this time a biholomorphism from $\Pi$ to $\Pi_D$ we have that
\begin{multline*}
    \nabla_x G_\Pi(x,y) \cdot \nabla^\bot \tgm_\Pi(x) = |T'(x)|^2 \nabla_x G_{\Pi_D}(T(x),T(y)) \cdot \nabla^\bot\tgm_{\Pi_D}(T(x))  \\ + \frac{\nabla_xG_{\Pi_D}(T(x),T(y))\cdot \psi(x)}{|T'(x)|^2},
\end{multline*}
where $\psi$ is another bounded function on $\Pi$ that has the same expression in terms of the conformal map $T$.

We specify now some properties of $\gamma_U$.
\begin{lemme}\label{lemmelimitesDeGamma}
We have that for any $x_0\in\partial U$, $\gamma_U(x,y) \tend{x,y}{x_0} + \infty$.
\end{lemme}
\begin{proof}
 Let $T : U \rightarrow D$ a biholomorphism. Relations \eqref{transpgamma} and \eqref{eqGammaD} yield that
\begin{equation}\label{DevGammaByT}
    \gamma_U(x,y) = -\frac{1}{2\pi}\ln\big(|T(x) - T(y)^*||T(y)\big) + \frac{1}{2\pi}\ln\frac{|T(x)-T(y)|}{|x-y|}.
\end{equation}
Obviously $|T(x) - T(y)^*||T(y)|\to |T(x_0) - T(x_0)^*||T(x_0)|=0$ when $x$ and $y$ go to $x_0$ so $-\frac{1}{2\pi}\ln\big(|T(x) - T(y)^*||T(y)|\big)$ goes to $+\infty$. Therefore we only need to obtain a lower bound for the other term.
By Theorem \ref{theoKW} the map $x,y \mapsto \frac{|T(x)-T(y)|}{|x-y|}$ is continuous and non zero on $\overline{U}\times \overline{U}$. Therefore, there exists a constant $C$ such that 
\begin{equation}\label{eqMajInfT}
    \frac{|T(x)-T(y)|}{|x-y|} > C,
\end{equation} 
for all $x,y\in\overline{U}$ and thus the lemma is proved.
\end{proof}
Noticing that for any exterior domain $\Pi$ and biholomorphism $T : \Pi\rightarrow \Pi_D$ there exists a neighborhood of $\partial \Pi \times \partial\Pi $ in $\Pi\times\Pi$ and a constant $C$ such that in this neighborhood relation \eqref{eqMajInfT} holds, the proof of the previous lemma holds for exterior domains too and thus,
\begin{equation}\label{limiteGammaPi}
\forall x_0\in\partial \Pi, \sp \gamma_\Pi(x,y) \tend{x,y}{x_0} + \infty.
\end{equation}
The following lemma gives explicit estimates of $d(x,\partial U)$, $d(y,\partial U)$ and $|x-y|$ when $\gamma_U(x,y) \rightarrow + \infty$.

\begin{lemme}\label{LemmeGammaMinoré} Let $k>0$ and $M \ge 0$ be constants. Let $\eps > 0 $. Assume that $(x,y)\in U\times U$ are such that  $$\frac{k}{2\pi}|\ln\eps| - M \le \gamma_U(x,y).$$ Then there exists a constant $C=C(M,U)$ such that 
\begin{equation*}
    \begin{cases}
        |x-y| \le C \eps^{k} \\
        d(x,\partial U) \le C \eps^{k} \\
        d(y,\partial U) \le C \eps^{k}.
    \end{cases}
\end{equation*}
\end{lemme}
\begin{proof}

Using relation \eqref{DevGammaByT} we have
\begin{equation*}
          \frac{k}{2\pi}|\ln\eps| - M \le -\frac{1}{2\pi}\ln\big(|T(x) - T(y)^*||T(y)|\big) + \frac{1}{2\pi}\ln\frac{|T(x)-T(y)|}{|x-y|}.
\end{equation*}
Recalling relation \eqref{eqMajInfT}, there exists a constant $C$ such that
\begin{equation*}
          \frac{k}{2\pi}|\ln\eps| - M \le -\frac{1}{2\pi}\ln\big(|T(x) - T(y)^*||T(y)|\big) +C,
\end{equation*}
and thus
\begin{equation*}
        |T(x)\overline{T(y)}-1| \le \eps^{k} e^{2\pi(C+M)} .
\end{equation*}
Moreover we have that for every $(a,b) \in D$, $|a-b| \le |1-a\overline{b}|$. Indeed, one can check that
\begin{equation*}
    |1-a\overline{b}|^2 - |a-b|^2 = (1-|b|^2)(1-|a|^2) > 0.
\end{equation*}
Therefore 
\begin{equation*}
        |T(x)-T(y)| \le \eps^k e^{2\pi(C-M)}
\end{equation*}
and relation \eqref{eqMajInfT} gives that
\begin{equation*}
    |x-y| \le C \eps^k.
\end{equation*}
It also yields that
\begin{equation*}
    (1-|T(x)|^2)(1-|T(y)|^2) \le  |T(x)\overline{T(y)}-1|^2 \le C \eps^{2k}.
\end{equation*}
That means that either $1-|T(x)|^2 \le \sqrt{C} \eps^k$, or $1-|T(y)|^2 \le \sqrt{C} \eps^k$.
We can assume without loss of generality that $1-|T(x)|^2 \le \sqrt{C} \eps^k$. We infer that $ 1-|T(x)| \le \sqrt{C} \eps^k $
and by the properties of the map $T^{-1}$ given in Theorem \ref{theoKW}, we conclude that 
$d(x,\partial U) \le C' \eps^k.$
Since $|x-y| \le C \eps^k$, that means that  $d(x,\partial U) \le C \eps^k$ and $d(y,\partial U) \le C \eps^k$.
\end{proof}

This lemma also stands in the case of an exterior domain $\Pi$ as follows.
\begin{lemme}\label{LemmeGammaMinoréExterior} Let $k>0$ and $M \ge 0$ be constants. Let $\eps > 0 $. Let $\mathcal{U} \subset \Pi$ be a bounded domain. Assume that $(x,y)\in \mathcal{U}\times\mathcal{U}$ are such that  $$\frac{k}{2\pi}|\ln\eps| - M \le \gamma_\Pi(x,y).$$ Then there exists a constant $C=C(M,\mathcal{U},\Pi)$ such that 
\begin{equation*}
    \begin{cases}
        |x-y| \le C \eps^{k} \\
        d(x,\partial\Pi) \le C \eps^{k} \\
        d(y,\partial\Pi) \le C \eps^{k}.
    \end{cases}
\end{equation*}
\end{lemme}
\begin{proof}
We argue as in Lemma \ref{LemmeGammaMinoré}. Relation \eqref{GreenPiD} gives that if $T:\Pi \rightarrow \Pi_D$ is a biholomorphism, then
\begin{equation*}
    \gamma_\Pi(x,y) = -\frac{1}{2\pi}\ln\big(|T(x) - T(y)^*||T(y)|\big) + \frac{1}{2\pi}\ln\frac{|T(x)-T(y)|}{|x-y|}.
\end{equation*}
Relation \eqref{eqMajInfT} holds true on the set $\overline{\mathcal{U}}\times\overline{\mathcal{U}}$, so we still have that
\begin{equation*}
        |T(x)\overline{T(y)}-1| \le \eps^{k} e^{2\pi(C+M)} .
\end{equation*}

We notice now that for every $(a,b) \in \Pi_D$ we also have that $|a-b| \le |1-a\overline{b}|$ since $(1-|b|^2)(1-|a|^2) > 0$. Thus
\begin{equation*}
    |x-y| \le C \eps^k
\end{equation*}
and
\begin{equation*}
    (1-|T(x)|^2)(1-|T(y)|^2) \le  |T(x)\overline{T(y)}-1|^2 \le C \eps^{2k}.
\end{equation*}
We have either $|T(x)|^2-1 \le \sqrt{C} \eps^k$, or $|T(y)|^2-1 \le \sqrt{C} \eps^k$, and by the same argument, recalling Theorem \ref{theoComplement}, we have that 
\begin{equation*}
    |T(x)|^2-1 \le \sqrt{C} \eps^k \Longrightarrow d(x,\partial \Pi) \le C \eps^k.
\end{equation*}
This proves the lemma.
\end{proof}

The next lemma give the formula for the exterior normal vector to $\partial U$.

\begin{lemme}\label{lemmeTranspVecteurNormal}
Let $x\in\partial U$. If $n_U(x)$ is the exterior normal unit vector to $\partial U$ in $x$, and $n_D(T(x))$ the exterior normal unit vector to $\partial D$ in $T(x)$, then
\begin{equation}\label{eqTranspVecteurNormal}
    n_U(x) = \frac{1}{|T'(x)|}\begin{pmatrix}\partial_1T_1(x) n^1_D(T(x)) + \partial_1T_2(x)n^2_D(T(x)) \\ \partial_1T_1(x) n^2_D(T(x)) - \partial_1T_2(x)n^1_D(T(x)) \end{pmatrix}.
\end{equation}
\end{lemme}
\begin{proof}
Since we chose $U$ and $T$ as in Theorem \ref{theoKW}, the map: \begin{equation*}
    \Gamma :  \R \to \partial U, \quad \Gamma(\theta)= T^{-1}(e^{i\theta})
\end{equation*}
is well defined and smooth. Therefore, denoting $x=\Gamma(\theta)$ we have that $T(x) = e^{i\theta}$ and
\begin{equation*}
      \Gamma'(\theta) = ie^{i\theta}(T^{-1})'(e^{i\theta}) = \frac{ie^{i\theta}}{T'(T^{-1}(e^{i\theta}))} = \frac{ie^{i\theta}}{T'(x)}.
\end{equation*}
Naturally, the exterior normal unit vector $n_D(e^{i\theta})$ to $\partial D$ in $e^{i\theta}$ is itself $e^{i\theta}$. Since an holomorphic map preserves the orientation, $ -i\frac{\Gamma'(\theta)}{|\Gamma'(\theta)|}$ is the exterior normal unit vector to $\partial U$ in $T^{-1}(e^{i\theta})=x$. Therefore
\begin{equation*}
    n_U(x) = -i\frac{\Gamma'(\theta)}{|\Gamma'(\theta)|} = \frac{e^{i\theta}}{T'(x)} |T'(x)| =\frac{ n_D(T(x)) \overline{T'(x)}}{|T'(x)|}.
\end{equation*}
We then compute the product 
\begin{equation*}
    n_D \overline{T'} = (n_D^1 + in_D^2)(\partial_1 T_1 - i \partial_1 T_2) = n_D^1\partial_1 T_1 + n_D^2\partial_1 T_2 + i(n_D^2\partial_1 T_1 - n_D^1\partial_1 T_2 )
\end{equation*}
and the lemma is now proved.
\end{proof}
We  extend the map $n_D$ to the interior of $D$ by the natural formula $n_D(x) = x$.  We then extend the map $n_U$ to $U$ by formula \eqref{eqTranspVecteurNormal}. The following lemma holds true.

\begin{lemme}\label{lemmeTgmOrthoSC}
There exists a constant $C$ such that for every $x \in U$, we have that
\begin{equation*}
    |\nabla^\perp \tgm_U(x) \cdot n_U(x) | \le C
\end{equation*}
\end{lemme}
\begin{proof}
We compute the scalar product $\nabla^\bot \tgm_U \cdot n_U$ by using relations \eqref{eqNablaTgmDev} and \eqref{eqTranspVecteurNormal}. We obtain
\begin{multline*}
    \nabla^\bot \tgm_U \cdot n_U = \begin{pmatrix}\ds \partial_1 \tgm_D  \partial_1 T_2 - \partial_2 \tgm_D  \partial_1 T_1  + \frac{+\partial_1 T_1 \partial_1^2 T_2 - \partial_1 T_2 \partial_1^2 T_1}{2\pi|T'|^2} \vspace{2mm} \\ \ds \partial_1 \tgm_D \partial_1 T_1 + \partial_2 \tgm_D  \partial_1 T_2  + \frac{\partial_1 T_1 \partial_1^2 T_1 + \partial_1 T_2 \partial_1^2 T_2}{2\pi|T'|^2}\end{pmatrix} \\ \cdot \frac{1}{|T'|} \begin{pmatrix}\partial_1T_1 n^1_D + \partial_1T_2n^2_D \\ \partial_1T_1 n^2_D - \partial_1T_2n^1_D \end{pmatrix} .
    \end{multline*}

We write  $$\nabla^\bot \tgm_U \cdot n_U \equiv \frac{1}{|T'|}A+\frac{1}{2\pi|T'|^3}B$$ with
\begin{equation*}
    A =  \begin{pmatrix}  \partial_1 \tgm_D \partial_1 T_2 - \partial_2 \tgm_D \partial_1 T_1 \\  \partial_1 \tgm_D \partial_1 T_1 + \partial_2 \tgm_D \partial_1 T_2 \end{pmatrix}\cdot  \begin{pmatrix}\partial_1T_1 n^1_D + \partial_1T_2n^2_D \\ \partial_1T_1 n^2_D - \partial_1T_2n^1_D \end{pmatrix} 
\end{equation*}
and 
\begin{equation*}
    B = \begin{pmatrix}\ds \partial_1 T_1 \partial_1^2 T_2 - \partial_1 T_2 \partial_1^2 T_1 \\ \ds \partial_1 T_1 \partial_1^2 T_1 + \partial_1 T_2 \partial_1^2 T_2\end{pmatrix} \cdot  \begin{pmatrix}\partial_1T_1 n^1_D + \partial_1T_2n^2_D \\ \partial_1T_1 n^2_D - \partial_1T_2n^1_D \end{pmatrix}.
\end{equation*}
We have that
\begin{align*}
        A & = n_D^1 \bigg[ \partial_1 \tgm_D \big(\partial_1 T_1 \partial_1T_2 - \partial_1 T_2\partial_1T_1 \big)  - \partial_2 \tgm_D\big(\partial_1 T_1 \partial_1T_1 + \partial_1 T_2\partial_1T_2 \big)  \bigg]
     \\ & \sp \sp+ n_D^2 \bigg[ \partial_1 \tgm_D \big(\partial_1 T_2 \partial_1T_2 + \partial_1 T_1\partial_1T_1 \big) + \partial_2 \tgm_D\big(-\partial_1 T_2 \partial_1T_1 + \partial_1 T_1\partial_1T_2 \big)  \bigg]
     \\ & =  - n_D^1 \partial_2 \tgm_D|T'|^2 + n_D^2\partial_1 \tgm_D |T'|^2
     \\ & = |T'|^2 n_D \cdot \nabla^\perp \tgm_D.
\end{align*}
We know from relation \eqref{eqTgmD} that $\tgm_D$ is a radial function and thus $n_D \cdot \nabla^\perp \tgm_D = 0$. So $A = 0$. We now compute $B$.
\begin{align*}
    B & =  \begin{pmatrix}\ds \partial_1 T_1 \partial_1^2 T_2 - \partial_1 T_2 \partial_1^2 T_1 \\ \ds \partial_1 T_1 \partial_1^2 T_1 + \partial_1 T_2 \partial_1^2 T_2\end{pmatrix} \cdot  \begin{pmatrix}\partial_1T_1 n^1_D + \partial_1T_2n^2_D \\ \partial_1T_1 n^2_D - \partial_1T_2n^1_D \end{pmatrix} \\
    & =  n^1_D ( \partial_1^2 T_1 (-2\partial_1 T_2\partial_1 T_1) + \partial_1^2 T_2 ((\partial_1T_1)^2 - (\partial_1T_2)^2))\\
    & \sp \sp + n^2_D ( \partial_1^2 T_1((\partial_1T_1)^2 - (\partial_1T_2)^2)  + 2\partial_1^2 T_2 \partial_1 T_2\partial_1 T_1)\\
    & = 2\pi \psi(T)\cdot n_D,
\end{align*}
where the map $\psi$ is defined by relation \eqref{DefPsi}, and is bounded. Since there exists a constant $m$ such that $|T'(x)| > m > 0$, we infer that there exists a constant $C$ such that
\begin{equation*}
    |\nabla^\perp \tgm_U(x) \cdot n_U(x) | \le C.
\end{equation*}
\end{proof}

Similarly, the normal vector $n_{\Pi_D}$ can be extended to $\Pi_D$ as a smooth function by the formula $n_{\Pi_D}(x) = -x$ for all $x\in\Pi_D$. We can then reproduce  Lemma \ref{lemmeTranspVecteurNormal} to extend $n_\Pi$ to the interior of $\Pi$ by the formula
\begin{equation*}
    n_\Pi(x) = \frac{1}{|T'(x)|}\begin{pmatrix}\partial_1T_1(x) n^1_{\Pi_D}(T(x)) + \partial_1T_2(x)n^2_{\Pi_D}(T(x)) \\ \partial_1T_1(x) n^2_{\Pi_D}(T(x)) - \partial_1T_2(x)n^1_{\Pi_D}(T(x)) \end{pmatrix}.
\end{equation*}

Lemma \ref{lemmeTgmOrthoSC} can be adapted to the exterior domain case in a straightforward manner. We obtain that for any exterior domain $\Pi$, and any bounded subset $\mathcal{U}\subset\Pi$, there exists a constant $C$ such that $\forall x \in \mathcal{U}$,
\begin{equation}\label{eqTgmOrthoPi}
    |\nabla^\perp \tgm_\Pi(x) \cdot n_\Pi(x) | \le C.
\end{equation}

\subsection{Inequalities for simply connected bounded domains and exterior domains}

We start with an estimate on the gradient of the Green's function. Let $\mathcal{U}$ be a bounded domain with $C^{2,\alpha}$ boundary. There exists a constant $C$ depending only on $\mathcal{U}$ such that
\begin{equation}\label{majDerGreen}
    \forall (x,y) \in \mathcal{U}\times\mathcal{U}, \sp x \neq y, \sp\sp |\nabla_xG_\mathcal{U}(x,y)| \le \frac{C}{|x-y|}.
\end{equation}
This estimate can be found in \cite{lichtenstein1919neuere}, see also \cite[Proposition 6.1]{iftimie_weak_2020}.

Now we can state the required inequalities for simply connected domains. We recall that $U$ is a simply connected bounded domain with $C^{2,\alpha}$ boundary.
\begin{lemme}\label{lemmeMajScalSimpleConnexe}
The following inequalities hold true for any $\kappa <1$:
\begin{equation*}
     \iint_{U\times U}\frac{1}{d(x,\partial U)^\kappa}|\nabla_xG_U(x,y) \cdot \nabla^\bot\tgm_U(x)| \dd x\dd y < \infty
\end{equation*}
and
\begin{equation*}
     \iint_{U\times U}\frac{1}{|x-y|^\kappa}|\nabla_xG_U(x,y) \cdot \nabla^\bot\tgm_U(x)| \dd x\dd y < \infty.
\end{equation*}
\end{lemme}

\begin{proof}
We start by denoting either $p(x,y) =\frac{1}{|x-y|^\kappa}$ or $p(x,y) = \frac{1}{d(x,\partial U)^\kappa}$. We use Lemma \ref{lemmeTransp}. Since $\psi$ is bounded, and since from Theorem \ref{theoKW} there exist constants $m$ and $M$ such that for every $x \in U$, $0<m<|T'(x)|<M$, we obtain that
\begin{multline*}
     \iint_{U\times U}p(x,y)|\nabla_xG_U(x,y) \cdot \nabla^\bot\tgm_U(x)| \dd x\dd y \\ \le  C\iint_{U\times U}p(x,y) [ |\nabla_x G_D(T(x),T(y)) \cdot \nabla^\bot\tgm_D(T(x))| + |\nabla_xG_D(T(x),T(y))| ]\dd x \dd y.
\end{multline*}
We now change variables using that $0<m<|T'(x)|<M$, to obtain that
\begin{multline*}
     \iint_{U\times U}p(x,y)|\nabla_xG_U(x,y) \cdot \nabla^\bot\tgm_U(x)| \dd x\dd y \\ \le C\iint_{D \times D} p(T^{-1}(x),T^{-1}(y))[|\nabla_xG_{D}(x,y) \cdot \nabla^\bot\tgm_{D}(x)| + |\nabla_x G_D(x,y)| ]\dd x \dd y.
\end{multline*}

Assume now that $ p(x,y) = \frac{1}{d(x,\partial U)^\kappa}$. By the properties of the map $T$, there exists a constant $C>0$ such that $d(T^{-1}(x),\partial U) \ge Cd(x,\partial D)$. We thus have that 
\begin{equation*}
    p(T^{-1}(x),T^{-1}(y)) = \frac{C}{d(T^{-1}(x),\partial U)^\kappa} \le \frac{C}{d(x,\partial D)^\kappa}.
\end{equation*}
In the light of this inequality, using relation \eqref{majDerGreen} and Lemma \ref{lemmeIneqIntermediaires} on $\mathcal{U} = D$ yields that 
\begin{equation*}
    \iint_{D \times D} p(T^{-1}(x),T^{-1}(y))|\nabla_x G_D(x,y)|\dd x \dd y \le C\iint_{D \times D} \frac{1}{d(x,\partial D)^\kappa|x-y|}<  \infty.
\end{equation*}
Recalling relation \eqref{eqMarchioro1} implies
\begin{equation*}
     \iint_{U\times U}\frac{1}{d(x,\partial U)^\kappa}|\nabla_xG_U(x,y) \cdot \nabla^\bot\tgm_U(x)| \dd x\dd y < \infty.
\end{equation*}

We assume next that  $p(x,y) = \frac{1}{|x-y|^\kappa}$. By the properties of the map $T$, we have for every $(x,y) \in D\times D$, $x \neq y$ that $|T^{-1}(x)-T^{-1}(y)|>C|x-y|$. Thus
\begin{equation*}
    p(T^{-1}(x),T^{-1}(y)) \le C p(x,y).
\end{equation*}
Using once again \eqref{majDerGreen} and Lemma \ref{lemmeIneqIntermediaires}, we have that 
\begin{equation*}
    \iint_{D \times D} p(T^{-1}(x),T^{-1}(y))|\nabla_x G_D(x,y)|\dd x \dd y \le C\iint_{D \times D} \frac{1}{|x-y|^{1+\kappa}}<  \infty,
\end{equation*}
and recalling relation \eqref{eqMarchioro2} we conclude that
\begin{equation*}
     \iint_{U\times U}\frac{1}{|x-y|^\kappa}|\nabla_xG_U(x,y) \cdot \nabla^\bot\tgm_U(x)| \dd x\dd y < \infty.
\end{equation*}
The lemma is proved.
\end{proof}

We now show these inequalities for $\Pi_D$. However the integral must be taken on a bounded subset.

\begin{lemme} \label{LemmeMajScalPiD}
Let $\mathcal{U}$ be any bounded subset of $\Pi_D$.
The following inequalities hold true for any $\kappa <1$:
\begin{equation*}
     \iint_{\mathcal{U}\times \mathcal{U}}\frac{1}{d(x,\partial \Pi_D)^\kappa}|\nabla_xG_{\Pi_D}(x,y) \cdot \nabla^\bot\tgm_{\Pi_D}(x)| \dd x\dd y < \infty
\end{equation*}
and
\begin{equation*}
     \iint_{\mathcal{U}\times \mathcal{U}}\frac{1}{|x-y|^\kappa}|\nabla_xG_{\Pi_D}(x,y) \cdot \nabla^\bot\tgm_{\Pi_D}(x)| \dd x\dd y < \infty.
\end{equation*}
\end{lemme}
\begin{proof}
We recall that $G_D$ and $G_{\Pi_D}$ have the same explicit expression, see relations \eqref{GreenD} and \eqref{GreenPiD}. Noticing that for $(x,y) \in \C\times\C^*$, $|x-y^*||y| = |x\overline y - 1|$ we obtain that for every $(x,y) \in \Pi_D \times \Pi_D$, $x\neq y$,
\begin{equation*}
     G_{D}\left(\frac{1}{x},\frac{1}{y}\right)= \frac{1}{2\pi} \ln \frac{\left|\frac{1}{x} - \frac{1}{y}\right|}{|\frac{1}{x\overline y} - 1 |} = \frac{1}{2\pi} \ln \frac{|x-y|}{|x\overline y-1|} = G_{\Pi_D}(x,y).
\end{equation*}
We can reproduce the proof of Lemma \ref{lemmeTransp} where we replace the biholomorphism  $T:U\to D$ with the map $T : \Pi_D \rightarrow D,$  $T(z) = 1/z$. This map is holomorphic and satisfies that $0<m<|T'(x)|<M$ and $|T''(x)|<M$ on $\mathcal{U}$. The calculations given in the proof of Lemma \ref{lemmeTransp} allow to obtain the following bound 
\begin{equation*}
    |\nabla_x G_{\Pi_D}(x,y) \cdot \nabla^\bot \tgm_{\Pi_D}(x)| \le C \left|\nabla_x  G_{D}\left(\frac{1}{x},\frac{1}{y}\right) \cdot \nabla^\bot\tgm_D\left(\frac{1}{x}\right)\right| + C  \left|\nabla_xG_{D}\left(\frac{1}{x},\frac{1}{y}\right)\right|.
\end{equation*}
for all $x,y\in\mathcal{U}$.

Thus for any $p : \C^2 \rightarrow \R$, we have that
\begin{multline*}
    \iint_{\mathcal{U}\times\mathcal{U}}p(x,y)|\nabla_xG_{\Pi_D}(x,y) \cdot \nabla^\bot\tgm_{\Pi_D}(x)| \dd x\dd y  \\ \le C\iint_{\mathcal{U}\times\mathcal{U}}p(x,y) \left(\left|\nabla_x  G_{D}\left(\frac{1}{x},\frac{1}{y}\right) \cdot \nabla^\bot\tgm_D\left(\frac{1}{x}\right)\right| +   \left|\nabla_xG_{D}\left(\frac{1}{x},\frac{1}{y}\right)\right|\right) \dd x \dd y.
\end{multline*}
Changing variables we obtain
\begin{multline*}
     \iint_{\mathcal{U}\times\mathcal{U}}p(x,y)|\nabla_xG_{\Pi_D}(x,y) \cdot \nabla^\bot\tgm_{\Pi_D}(x)| \dd x\dd y \\ \le C\iint_{T(\mathcal{U}) \times T(\mathcal{U})} p\left(\frac{1}{x},\frac{1}{y}\right)\left(|\nabla_xG_{D}(x,y) \cdot \nabla^\bot\tgm_{D}(x)| + |\nabla_x G_{D}(x,y)| \right) \dd x \dd y.
\end{multline*}
The end of the proof is very similar to Lemma \ref{lemmeMajScalSimpleConnexe}. We start by showing that in the case  $p(x,y) = \frac{1}{d(x,\partial \Pi_D)^\kappa} = \frac{1}{(|x|-1)^\kappa}$, we have for every $x,y \in T(\mathcal{U})\times T(\mathcal{U})$ that 
\begin{equation*}
    p\left(\frac{1}{x},\frac{1}{y}\right) =   \frac{|x|^\kappa}{(1-|x|)^\kappa} \le \frac{1}{d(x,\partial D)^\kappa}.
\end{equation*}

We can use relation \eqref{majDerGreen} and Lemma \ref{lemmeIneqIntermediaires} on $ T(\mathcal{U})$ which is a bounded domain, to observe that
\begin{equation*}
    \iint_{T(\mathcal{U}) \times T(\mathcal{U})} \frac{1}{d(x,\partial D)^\kappa} |\nabla_x G_{D}(x,y)| \dd x \dd y < \infty.
\end{equation*}

Relation \eqref{eqMarchioro1} implies that
\begin{equation*}
   \iint_{T(\mathcal{U}) \times T(\mathcal{U})} \frac{1}{d(x,\partial D)^\kappa}|\nabla_xG_{D}(x,y) \cdot \nabla^\bot\tgm_{D}(x)| \dd x \dd y < \infty.
\end{equation*}
We proved that
\begin{equation*}
    \iint_{\mathcal{U}\times\mathcal{U}}\frac{1}{d(x,\partial \Pi_D)^\kappa}|\nabla_xG_{\Pi_D}(x,y) \cdot \nabla^\bot\tgm_{\Pi_D}(x)| \dd x\dd y < \infty.
\end{equation*}

Similarly, in the case $p(x,y) = \frac{1}{|x-y|^\kappa}$, we have that 
\begin{equation*}
    p\left(\frac{1}{x},\frac{1}{y}\right) = \frac{|xy|^\kappa}{|y-x|^\kappa} \le \frac{1}{|y-x|^\kappa}.
\end{equation*}
Since $\mathcal{U}$ is bounded, Lemma \ref{lemmeIneqIntermediaires} apply, and using relation \eqref{majDerGreen} yields that
\begin{equation*}
    \iint_{\mathcal{U}\times\mathcal{U}}\frac{1}{|y-x|^\kappa}|\nabla_xG_{D}(x,y)| \dd x\dd y < \infty.
\end{equation*}
Relation \eqref{eqMarchioro2} implies that
\begin{equation*}
   \iint_{T(\mathcal{U}) \times T(\mathcal{U})} \frac{1}{|x-y|^\kappa}|\nabla_xG_{D}(x,y) \cdot \nabla^\bot\tgm_{D}(x)| \dd x \dd y < \infty.
\end{equation*}
We conclude that
\begin{equation*}
    \iint_{\mathcal{U}\times\mathcal{U}} \frac{1}{|x-y|^\kappa}|\nabla_xG_{\Pi_D}(x,y) \cdot \nabla^\bot\tgm_{\Pi_D}(x)| \dd x\dd y < \infty.
\end{equation*}
The lemma is now proved.
\end{proof}

We conclude this section with the same inequalities for any exterior domain $\Pi$.
\begin{lemme}\label{lemmeMajScalExterior} Let $\mathcal{U}$ be any bounded subset of $\Pi$.
The following inequalities hold true for any $\kappa <1$:
\begin{equation*}
     \iint_{\mathcal{U}\times \mathcal{U}}\frac{1}{d(x,\partial \Pi)^\kappa}|\nabla_xG_{\Pi}(x,y) \cdot \nabla^\bot\tgm_{\Pi}(x)| \dd x\dd y < \infty
\end{equation*}
and
\begin{equation*}
     \iint_{\mathcal{U}\times \mathcal{U}}\frac{1}{|x-y|^\kappa}|\nabla_xG_{\Pi}(x,y) \cdot \nabla^\bot\tgm_{\Pi}(x)| \dd x\dd y < \infty.
\end{equation*}
\end{lemme}
\begin{proof}
The proof follows the same outline as the proofs of Lemmas \ref{lemmeMajScalSimpleConnexe} and \ref{LemmeMajScalPiD}.
\end{proof}

\section{Multiply connected domain}\label{sectionMultiConnected}
We work now with the multiply connected domain $\Omega$.
There exists a compact set $K$ such that $\Omega\setminus K$ has exactly $n+1$ connected components $V_0,\ldots,V_n$ that satisfy $d(V_i,\Gamma_j) > 0$ for every $i\neq j$. For example one can take $V_i$ to be the $\eps$-neighborhood of $\Gamma_i$, for $\eps$ small enough. Thus we also have that
\begin{equation}\label{decompOmegaVj}
    \Omega = K \cup \left( \bigcup_{j=0}^m V_j\right).
\end{equation}
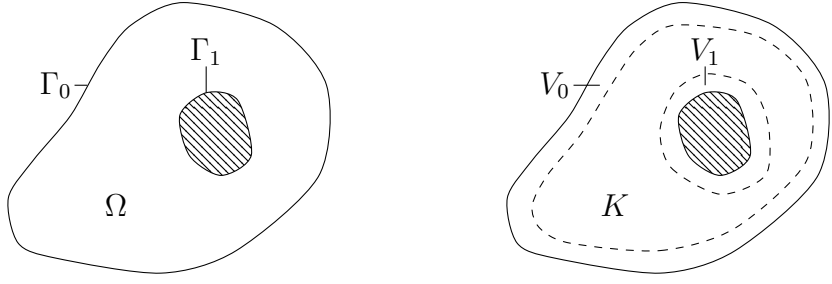
\begin{figure}\centering
    \begin{tikzpicture}
    \draw plot [smooth cycle]coordinates  {(1,1) (2.4,0.8) (3.4, 1.1) (4.5, 2.1) (4.6, 3.3) (3.9,4.1) (2.7,4.4) (2.0, 4.1) (1.6,3.6) (1.2, 2.9) (0.7, 2.3) (0.37,1.8) (0.5,1.2)};
    \draw[pattern=north west lines]  plot [smooth cycle]coordinates  {(2.65,2.90) (2.78,2.36) (3.2,2.1) (3.6,2.37) (3.42,3.10) (3,3.2)}node at (1.8,1.7) {$\Omega$} node at (1,3.3) {$\Gamma_0$} node at (3,3.75)  {$\Gamma_1$};
    \draw (3,3.55) -- (3,3.2);
    \draw(1.25,3.3) --(1.42,3.3);
    \end{tikzpicture}
    \hspace{2cm}
    \begin{tikzpicture}
    \draw plot [smooth cycle]coordinates  {(1,1) (2.4,0.8) (3.4, 1.1) (4.5, 2.1) (4.6, 3.3) (3.9,4.1) (2.7,4.4) (2.0, 4.1) (1.6,3.6) (1.2, 2.9) (0.7, 2.3) (0.37,1.8) (0.5,1.2)};
    \draw[dashed] plot [smooth cycle]coordinates  {(1.2,1.2) (2.4,1.1) (3.4, 1.4) (4.35, 2.3) (4.35, 3.35) (3.7,3.95) (2.6,4.1) (2.1, 3.8) (1.75,3.3) (1.35, 2.65) (0.9, 2.1) (0.7,1.7) (0.8,1.4)} node at (1.8,1.7) {$K$};
    \draw[dashed]  plot [smooth cycle]coordinates  {(2.4,2.70) (2.55,2.2) (3.2,1.85) (3.8,2.1) (3.85,2.7) (3.72,3.12) (3.5,3.35) (3,3.45) (2.6,3.22)};
    \draw[pattern=north west lines]  plot [smooth cycle]coordinates  {(2.65,2.90) (2.78,2.36) (3.2,2.1) (3.6,2.37) (3.42,3.10) (3,3.2)};
    \draw node at (3,3.75) {$V_1$};
    \draw (3,3.55) -- (3,3.3);
    \draw node at (1,3.3) {$V_0$};
    \draw(1.25,3.3) --(1.6,3.3);
    \end{tikzpicture}
    \caption{Decomposition \eqref{decompOmegaVj}.}
    \label{fig2}
\end{figure}

\subsection{Biot-Savart law}
Let $\omega$ be a fixed function on $\Omega$. Obtaining the velocity $u$ in terms of the vorticity $\omega$ is solving the following problem
\begin{equation}\label{inverseProblem}
    \begin{cases}
    \curl u = \omega, & \text{ in }\Omega\\
    \nabla \cdot u = 0, & \text{ in }\Omega\\
    u \cdot n = 0, & \text{ on }\partial\Omega.
    \end{cases}
\end{equation}
As in the simply connected case, a particular solution is given by
\begin{equation*}
    u(x) = \int_\Omega \nabla_x^\bot G_\Omega(x,y)\omega(y)\dd y.
\end{equation*}
Since \eqref{inverseProblem} is linear, the general solution is given by this particular solution plus the general solution of the homogeneous problem. The solution is of the form (see \cite{iftimie_weak_2020})
\begin{equation*}
    u(x,t) = \int \nabla_x^\bot G_\Omega(x,y)\omega(y,t)\dd y + \sum_{j=1}^m c_{j,\omega}(t) \nabla^\bot w_j(x)
\end{equation*}
where 
\begin{equation*}
c_{j,\omega}(t)=\int w_j(x)\omega(x,t)\dd x+\xi_j,
\end{equation*}
$\xi_j$ is the circulation of the velocity $u$ on $\Gamma_j$ and $w_j : \Omega \mapsto \R$ are the harmonic measures defined by 
\begin{equation*}
    \begin{cases}
    \Delta w_j = 0 & \text{ in }\Omega\\
    w_j = \delta_{j,l}  & \text{ on }\Gamma_l, 0 \leq l \leq n.
    \end{cases}
\end{equation*}

The vector fields
\begin{equation*}
    \beta_j(x) = \nabla^\bot w_j(x)
\end{equation*}
are called harmonic vector fields.

In the case of a discrete vorticity
\begin{equation*}
\omega(t)=\sum_{j=1}^N a_j\delta_{x_j(t)}
\end{equation*}
we define the point vortex dynamics in multiply connected bounded domains as follows:
\begin{equation*}   \forall 1 \le i \le N, \sp \der{x_i(t)}{t} = \sum_{\substack{j= 1\\ j\neq i}}^N \nabla_x^\bot G_\Omega(x_i(t),x_j(t))a_j + \nabla_x^\bot\gamma_\Omega(x_i(t),x_i(t)) a_i + \sum_{j=1}^m c_j(t) \beta_j(x)
\end{equation*}
where
\begin{equation*}
c_j(t)=\xi_j+\sum_{k=1}^N a_k w_j(x_k(t)).
\end{equation*}

Let us observe that by the Kelvin theorem, the circulations $\xi_j$ are constant in time. They are therefore prescribed at the initial time. The harmonic measures $w_j$ being smooth, we observe that the functions   $c_j : \R^+ \rightarrow \R$  are bounded.

\subsection{Inequalities for multiply connected domains}
We know from \cite[Proposition 6.1]{iftimie_weak_2020} that  the inequality
\begin{equation}\label{majGreen}
    |G_\Omega(x,y)| \le C(1 + |\ln|x-y||)
\end{equation}
holds true for bounded domains. We also recall relation \eqref{majDerGreen}:
\begin{equation*}
    |\nabla_xG_\Omega(x,y)| \le \frac{C}{|x-y|}.
\end{equation*}
We combine this with Lemma \ref{lemmeIneqIntermediaires}. Since $\Omega$ is bounded, for any $\kappa <1$ we have that
\begin{equation}\label{majGreenPuissance}
   \iint_{\Omega\times\Omega} \frac{1}{|x-y|^\kappa} |\nabla_xG_\Omega(x,y)| \dd x \dd y < \infty
\end{equation}
and
\begin{equation}\label{majGreenPuissanceDistance}
   \iint_{\Omega\times\Omega} \frac{1}{d(x,\partial\Omega)^\kappa} |\nabla_xG_\Omega(x,y)| \dd x \dd y < \infty.
\end{equation}
The following proposition gives an estimate of the map $\tgm_{\Omega}$ near the boundary.
\begin{prop}[Gustafsson \cite{gustafsson1979motion} Proposition 3.3]\label{PropEquivTGM}
Denoting by $K_j$ the connected components of $\R^2 \setminus \Omega$ and $d_j(x) = \inf \{|x-y|, y \in K_j \}$, $D_j(x) = \sup \{|x-y|, y \in K_j \}$, we have that 
\begin{equation*}
    \forall x \in \Omega, \sp \sp \ln d(x,\partial\Omega) \le -2\pi\tgm_\Omega(x) \le \min_j \ln \frac{4d_j(x)}{1-\frac{d_j(x)}{D_j(x)}}. 
\end{equation*}
\end{prop}

Clearly, $d_j(x) < D_j(x)$ for every $x \in \Omega$, so by compactness
\begin{equation*}
    \sup \left\{ \frac{d_j(x)}{D_j(x)}, x \in \Omega \right\} < 1.
\end{equation*}
That means that there exists a constant $C_1$ depending only on $\Omega$ such that for every $x\in \Omega$
\begin{equation*}
    \ln d(x,\partial\Omega) \le -2\pi\tgm_\Omega(x) \le \min_j \ln d_j + C_1
\end{equation*}
and thus
\begin{equation}\label{EquivTGM}
    \ln d(x,\partial\Omega) \le -2\pi\tgm_\Omega(x) \le \ln d(x,\partial\Omega) + C_1.
\end{equation}
In particular, 
\begin{equation*}
    \tgm_\Omega(x)  \sim -\frac{1}{2\pi} \ln d(x,\partial\Omega)\quad \text{as }x\rightarrow \partial\Omega
\end{equation*}
and 
\begin{equation}\label{tgmMinoree}
    \inf_\Omega \tgm_\Omega = \min_\Omega \tgm_\Omega > - \infty.
\end{equation}
In addition, we state in the following proposition an estimate of $\nabla \tgm_\Omega$ near the boundary.
\begin{prop}[Gustafsson \cite{gustafsson1979motion} Proposition 3.5]
There exists a constant $C$ such that for every $x\in\Omega$,
\begin{equation}\label{MajDerTGM}
    |\nabla \tgm_\Omega(x)| \le \frac{C}{d(x,\partial\Omega)}
\end{equation}
Moreover, we can take $C = \frac{1}{2\pi}$ if $\Omega$ is simply connected.
\end{prop}

We now want to compare the map $\gamma_\Omega$ near the boundary $\Gamma_j$, to the map $\gamma_{\Omega_j}$.
\begin{lemme}\label{LemmeDiffGammabounded}
For any $0\le j \le m$, the map $\gamma_\Omega - \gamma_{\Omega_j}$ is bounded in $V_j\times\Omega$.
\end{lemme}
\begin{proof}
We have that
\begin{equation*}
    \gamma_\Omega(x,y) - \gamma_{\Omega_j}(x,y) =  G_\Omega(x,y) - G_{\Omega_j}(x,y)
\end{equation*}

We fix $x\in V_j$, and we define $F(y) = G_\Omega(x,y) - G_{\Omega_j}(x,y)$. It satisfies that

\begin{equation*}
    \begin{cases}
    \Delta_y F(y) = 0, & \text{ on $\Omega$,} \\
    F(y)= 0, & \text{ on ${\Gamma_j}$,} \\
    |F(y)| \leq C& \text{ on ${\Gamma_k}$, $k\neq j$,}
    \end{cases}
\end{equation*}
where 
\begin{equation*}
    C = \sup_{\substack{x \in V_j \\ k\neq j \\y \in \Gamma_k}} |G_{\Omega_j}(x,y)| < \infty
\end{equation*} is a constant that does not depend on $x\in V_j$. The supremum is finite since $d(V_j,\Gamma_k) > 0$ for each $k\neq j$. Therefore by the maximum principle, \begin{equation*}
    |\gamma_\Omega(x,y) - \gamma_{\Omega_j}(x,y)| = |F(y)| \le \max_{y\in \partial\Omega} |F(y)| \le C
\end{equation*}
for every $(x,y)\in V_j\times\Omega$.
\end{proof}
Observe that, for $x\in V_j$, the map $\widetilde F(y) = \nabla_x\gamma_\Omega(x,y) - \nabla_x \gamma_{\Omega_j}(x,y) $ satisfies the exact same problem: its Laplacian vanishes over $\Omega$,  $\widetilde F(y)= 0$ on $\Gamma_j$ and it is bounded on $\Gamma_k$ by a map $C(x)$ that is itself bounded in $V_j$. Hence the map $\nabla_x\gamma_\Omega - \nabla_x\gamma_{\Omega_j}$ is also bounded in $V_j\times\Omega$. Since $V_j\times V_j \subset V_j\times \Omega$, we can set $y=x$ in those inequalities and obtain similar bounds for $\tgm_\Omega - \tgm_{\Omega_j}$ and $\nabla \tgm_\Omega - \nabla \tgm_{\Omega_j}$. In this manner we obtain the following corollary.
\begin{coro}\label{LemmeDiffTGMbounded}
For any $0\le j \le m$, we have that 
\begin{itemize}
    \item the map $\nabla_x\gamma_\Omega - \nabla_x\gamma_{\Omega_j}$ is bounded in $V_j\times\Omega$.
    \item the map $\tgm_\Omega - \tgm_{\Omega_j}$ is bounded in $V_j$.
    \item the map $\nabla\tgm_\Omega - \nabla\tgm_{\Omega_j}$ is bounded in $V_j$.
\end{itemize}
\end{coro}
Lemma \ref{lemmelimitesDeGamma} is stated for bounded simply-connected domains, it is also true for exterior domains, see relation \eqref{limiteGammaPi}.
Lemma \ref{LemmeDiffGammabounded} allows to prove it for multiply connected domains.
\begin{coro}\label{coroLimiteGamma}
We have that for any $x_0\in\partial\Omega$, $\gamma_\Omega(x,y) \tend{x,y}{x_0\in \partial\Omega} + \infty$.
\end{coro}
We now combine Lemma \ref{lemmeTgmOrthoSC} and Lemma \ref{LemmeDiffTGMbounded} to obtain the following result.
\begin{coro}\label{coroTgmOrtho}
There exists a constant  $C$ such that for every $x \in \Omega$ and every $1 \le j\le m$,
\begin{equation*}
    |\nabla\tgm_\Omega(x) \cdot \beta_j(x)| \le C.
\end{equation*}
\end{coro}
\begin{proof}
Let $1\le j\le m$ and $0 \le k \le m$. In a neighborhood of $\Gamma_k$, we decompose $$\beta_j(x) \equiv \beta^1_j(x) n_{\Omega_k}(x) + \beta^2_j(x) n^\perp_{\Omega_k}(x).$$
Since $\beta_j$ is tangent to the boundary, $\beta^1_j(x) = 0$ when $x \in \Gamma_k$. Since $\beta_j$ is smooth, there exists a constant $C_{j,k}$ such that $|\beta^1_j(x)| \le C_{j,k} d(x,\Gamma_k)$ in a neighborhood of $\Gamma_k$. Recalling relation \eqref{MajDerTGM} and provided that the neighborhood is sufficiently small so that  $d(x,\partial\Omega) = d(x,\Gamma_k)$, we have that
\begin{equation*}
    |\beta^1_j(x)\nabla\tgm_\Omega(x) \cdot  n_{\Omega_k}(x)| \le C_{j,k} d(x,\Gamma_k) \frac{C}{d(x,\partial\Omega)} \le C_{j,k}.
\end{equation*}
If $k=0$, we apply Lemma \ref{lemmeTgmOrthoSC}. If $k\neq 0$, we use relation \eqref{eqTgmOrthoPi}. In both cases, it yields that
\begin{equation*}
    |\nabla^\perp\tgm_{\Omega_k}(x) \cdot n_{\Omega_k}(x))| \le C_k
\end{equation*}
in a neighborhood of $\Gamma_k$. Thus by Corollary \ref{LemmeDiffTGMbounded}, and since $n_{\Omega_k}$ is bounded in that neighborhood,
\begin{equation*}
    |\nabla^\perp\tgm_{\Omega}(x) \cdot n_{\Omega_k}(x))|  \le C_{k}.
\end{equation*}
Consequently, since $\beta_j$ is bounded,
\begin{equation*}
    |\beta^2_j(x) \nabla\tgm_{\Omega}(x) \cdot n^\perp_{\Omega_k}(x))|  \le C_{j,k}.
\end{equation*}
Therefore on this neighborhood of the boundary $\Gamma_k$, there exists a constant $C_{j,k}$ such that
\begin{equation*}
    |\nabla^\perp\tgm_\Omega(x) \cdot \beta_j(x)| \le C_{j,k}
\end{equation*}
Outside of each of these neighborhoods, we know that the maps $\nabla\tgm_{\Omega}$ and $\beta_j$ are bounded. Therefore, as there are a finite number of boundaries $\Gamma_k$, and of maps $\beta_j$, there exists a constant $C$ depending only on $\Omega$ such that for every $x \in \Omega$,
\begin{equation*}
    |\nabla^\perp\tgm_\Omega(x) \cdot \beta_j(x)| \le C.
\end{equation*}

\end{proof}

We can now extend Lemma \ref{lemmeMajScalSimpleConnexe} to the case of multiply connected domains.
\begin{lemme}\label{lemmeMajScalGenerale}
The following inequalities hold true for any $\kappa <1$:
\begin{equation*}
     \iint_{\Omega\times\Omega}\frac{1}{|x-y|^\kappa}|\nabla_xG_\Omega(x,y) \cdot \nabla^\bot\tgm_\Omega(x)| \dd x\dd y < \infty.
\end{equation*}
and
\begin{equation*}
     \iint_{\Omega\times\Omega}\frac{1}{d(x,\partial\Omega)^\kappa}|\nabla_xG_\Omega(x,y) \cdot \nabla^\bot\tgm_\Omega(x)| \dd x\dd y < \infty.
\end{equation*}
\end{lemme}
\begin{proof}
Let us introduce the map $h$ defined by
\begin{equation*}
    h(x,y) = \nabla_xG_\Omega(x,y) \cdot \nabla^\bot\tgm_\Omega(x).
\end{equation*}
We must show that $ph \in L^1(\Omega\times\Omega)$ for $p(x,y) = \frac{1}{|x-y|^\kappa}$ and also for $p(x,y) = \frac{1}{d(x,\partial\Omega)^\kappa}$. We split the integral using the decomposition \eqref{decompOmegaVj}, pictured in Figure \ref{fig2}.

First, there exists a constant $C$ such that $|\nabla^\bot\tgm_\Omega(x)| \le C$ on $K$. Relations \eqref{majGreenPuissance} and \eqref{majGreenPuissanceDistance} thus imply that $ph \in L^1(K\times \Omega)$ for both expression of $p$. Now we must prove that $p h \in L^1(V_j\times \Omega)$ for every $0\le j\le m$. We fix $0\le j \le m$. By Corollary \ref{LemmeDiffTGMbounded} as well as relations \eqref{majGreenPuissance} and \eqref{majGreenPuissanceDistance}, we know that proving $p h \in L^1(V_j\times \Omega)$ is equivalent to proving that $p h_1 \in L^1(V_j\times \Omega)$, with
\begin{equation*}
    h_1(x,y) = \nabla_xG_\Omega(x,y) \cdot \nabla^\bot\tgm_{\Omega_j}(x).
\end{equation*}
Let us introduce $$h_2(x,y) = \nabla_xG_{\Omega_j}(x,y)\cdot\nabla^\bot\tgm_{\Omega_j}(x).$$ We have that $h_1 - h_2 \in C^\infty(\Omega\times\Omega)$ and that $\Delta_y (h_1 - h_2) = 0$. The maximum principle yields
\begin{equation}\label{h1-h_2}
    \forall (x,y)\in V_j\times \Omega, \sp |h_1(x,y) - h_2(x,y)| \le \sup_{y \in \partial\Omega} |h_1(x,y) - h_2(x,y)|.
\end{equation}
However, $h_1(x,y) = 0$ when $y \in \partial\Omega$, since  $\forall x \in \Omega, G_\Omega(x,y) =0$ when $y\in\partial\Omega$. And similarly, $h_2(x,y) = 0$ when $y \in \partial\Omega_j$. Thus,
\begin{equation*}
    \sup_{y \in \partial\Omega} |h_1(x,y) - h_2(x,y)| = \sup_{\substack{k\neq j \\ y \in \Gamma_k}} |h_2(x,y)|.
\end{equation*}
We need to bound $h_2(x,y)$ when $x \in V_j$ and $y\in \Gamma_k$. We decompose
\begin{equation*}
    \nabla_xG_{\Omega_j}(x,y) \equiv g_1(x,y)n_{\Omega_j(x)} + g_2(x,y) n^\perp_{\Omega_j(x)}
\end{equation*}
where $n_{\Omega_j(x)}$ is defined in Section \ref{sectionIntermediateLemmas}. We have that $g_2(x,y) = 0$ when $x \in \Gamma_j$ since $G_{\Omega_j}(x,y) = 0$ for every $(x,y) \in \Gamma_j \times \Gamma_k$ so $\nabla_x G_{\Omega_j}(x,y)$ is normal to the boundary $\Gamma_j$. By Theorem \ref{lemmeGreenReguliere}, $G_{\Omega_j}$ is $C^2$ up to the boundary except on the diagonal, thus there exists a constant $C$ independent of $y \in \Gamma_k$ such that $|g_2(x,y)| \le C d(x,\Gamma_j)$ for all $x\in V_j$ and $y\in\Gamma_k$. Using relation \eqref{MajDerTGM}, we have that $$|g_2(x,y) n^\perp_{\Omega_j(x)} \cdot \nabla^\bot\tgm_{\Omega_j}(x)| \le C.$$ Using Lemma \ref{lemmeTgmOrthoSC} if $j=0$ and relation \eqref{eqTgmOrthoPi} if $j\neq 0$, we have that 
$$|g_1(x,y) n_{\Omega_j(x)} \cdot \nabla^\bot\tgm_{\Omega_j}(x)| \le C$$
for all $x\in V_j$ and $y\in\Gamma_k$. So there exists a constant independent of $x\in V_j$ and $y \in \Gamma_k$ such that   $|h_2(x,y)| \le C$.
Thus there exists a constant $C$ such that 
\begin{equation}\label{h1-h2bord}
    \sup_{y \in \partial\Omega} |h_1(x,y) - h_2(x,y)| \le C.
\end{equation}
Therefore, since $p \in L^1(\Omega\times\Omega)$, relations \eqref{h1-h_2} and \eqref{h1-h2bord} yield that $p(h_1 -h_2) \in L^1(V_j\times\Omega)$. 

If $j=0$, we now apply Lemma \ref{lemmeMajScalSimpleConnexe} to $\Omega_0$. If $j \neq 0$ we apply Lemma \ref{lemmeMajScalExterior} to the domain $\Omega_j$ and to its bounded subset $\Omega$. In both cases, we get that $p h_2$ is integrable on $\Omega\times\Omega$. Therefore $ph_1 \in L^1(V_j\times\Omega)$ and thus $ph\in L^1(V_j\times\Omega)$. This completes the proof of the lemma.
\end{proof}

\section{Completion of the proof of Theorem \ref{mainresult}}\label{sectionPreuve}

In this section we denote by $G$, $\gamma$ and $\tgm$ the maps associated to the domain $\Omega$ in order to lighten the notations as there should be no ambiguity. 

\subsection{Construction of a regularized dynamic}\label{sectionRegDynamic}

We need to construct a dynamics that is well defined for every time and which is the same as the point vortex dynamics  as long as no point vortex is close to the boundary and no two point vortices are close to each other. More precisely, we need to construct two functions $G_\eps$ and $\tgm_\eps$ such that the dynamics
\begin{equation}
\label{regularizedDynamic}    \forall 1 \le i \le N, \sp \der{x^\eps_i(t)}{t} = \sum_{\substack{j= 1\\ j\neq i}}^N \nabla_x^\bot G_\eps(x^\eps_i(t),x^\eps_j(t))a_j + \frac{1}{2}\nabla^\bot\tgm_\eps(x^\eps_i(t)) a_i +  \sum_{j=1}^m c_j(t) \beta_j(x^\eps_i(t))
\end{equation}
is well defined in $\overline{\Omega}$ for every time. It suffices that $G_\eps\in C^2(\overline{\Omega}\times\overline{\Omega})$ and $\tgm_\eps\in C^1(\overline{\Omega})$ and that $\nabla_x^\bot G_\eps$ and $\nabla^\bot\tgm_\eps$ are tangent to $\partial\Omega$ when the first variable is at the boundary. Moreover, we want to ensure that the following implication is true for every $(x,y) \in \Omega\times\Omega$,
\begin{equation}\label{RegAvantTauEps}
\left\{\begin{aligned}
        |G_{\R^2}(x,y)| &< \frac{1}{2\pi} |\ln \eps| 
        \\ |\tgm(x)| &< \frac{1}{2\pi} |\ln \eps | 
        \\ |\gamma(x,y)| &< \frac{1}{2\pi} |\ln \eps|
    \end{aligned}\right.
    \sp \sp \Longrightarrow \sp \sp   
\left\{\begin{aligned}
    G_\eps(x,y) &= G(x,y) \\
    \tgm_\eps(x) &= \tgm(x).
\end{aligned}\right.
\end{equation}
This ensures that the maps $G_\eps$, and $\tgm_\eps$ are good approximations of the maps $G$ and $\tgm$ when $\eps$ goes to 0.

In order to have a proper control over the regularized maps, we also want to ensure that
\begin{equation}\label{PropRegFunctions}
\left\{\begin{aligned}
    |\tgm_\eps(x)| &\le  |\tgm(x)| \\
    |\nabla \tgm_\eps(x)| &\le |\nabla \tgm(x)| \\
    |G_\eps(x,y)| &\le  |G(x,y)| \\
    |\nabla_x G_\eps(x,y)| &\le \frac{C}{|x-y|}
\end{aligned}\right.
\end{equation}
for a constant $C$ independent of $\eps$.

We consider $f_\eps \in C^\infty(\R,\R)$ an odd map such that
\begin{equation*}
\begin{cases}
    f_\eps(r) = r, & \forall |r| < \frac{1}{2\pi}|\ln\eps| \\
    f_\eps(r) = L_\eps, & \forall r > \frac{1}{2\pi} |\ln\eps| + 1 \\
    0 \le f_\eps'(r) \le 1 , &\forall r \in \R
\end{cases}
\end{equation*}
for some constant $L_\eps$.

\subsubsection*{Construction of $\tgm_\eps$}
We define the regularized Robin function as
\begin{equation*}
    \tgm_\eps(x) = f_\eps(\tgm(x))
\end{equation*}
so that
\begin{equation*}
    \nabla\tgm_\eps(x) = \nabla\tgm(x)f_\eps'(\tgm(x)).
\end{equation*}
We recall that according to Proposition \ref{PropEquivTGM}, $\tgm(x) \tend{x}{\partial\Omega} + \infty$ and thus $\tgm_\eps(x) \tend{x}{\partial\Omega} L_\eps $ and $\nabla\tgm_\eps(x) \tend{x}{\partial\Omega} 0$. Therefore $\tgm_\eps \in C^1(\overline{\Omega})$ and $\nabla^\perp \tgm_\eps$ is indeed tangent to the boundary since it vanishes at the boundary, and satisfies both
\begin{equation*}
    |\tgm_\eps(x)| \le  |\tgm(x)|
\end{equation*}
and
\begin{equation*}
     |\nabla \tgm_\eps(x)| \le |\nabla \tgm(x)|.
\end{equation*}

\subsubsection*{Construction of $G_\eps$}
We define the regularized Green's function for $(x,y) \in \overline{\Omega}\times\overline{\Omega}$ as follows:
\begin{equation*}
    \begin{cases} 
        G_\eps(x,y) = f_\eps(G_{\R^2}(x,y)) + f_\eps(\gamma(x,y)) & \text{ if } (x,y) \in \Omega\times\Omega, \sp x \neq y\\
        G_\eps(x,y) = 0 & \text{ if } x \in \partial\Omega \text{ or } y \in \partial\Omega \\
        G_\eps(x,x) = -L_\eps  + f_\eps(\tgm(x)) & \text{ if } x \in \Omega.
    \end{cases}
\end{equation*}
Let us notice straight away that $G_\eps(x,y)=G_\eps(y,x)$.

 We collect some properties of $G_\eps$ in the following lemma.
\begin{lemme}
We have that  $G_\eps \in C^2(\overline{\Omega}\times\overline{\Omega})$ and that $\nabla_x^\perp G_\eps(x,y)$ is tangent to $\partial\Omega$ when $x\in\partial\Omega$. Moreover,  $|G_\eps(x,y)|\leq |G(x,y)|$ and there exists a constant $C$ independent of $\eps$ such that  $|\nabla_x G_\eps(x,y)|\leq\frac{C}{|x-y|}$.
\end{lemme}
\begin{proof}
We start by proving that $G_\eps \in C^1(\overline{\Omega}\times\overline{\Omega})$. Since $f_\eps \in C^\infty(\R)$, the map $G_\eps$ is clearly $C^\infty$ on the set $\Omega \times \Omega \setminus \{ x = y \}$. 

We show first the continuity over $\overline{\Omega}\times\overline{\Omega}$. From relation \eqref{EquivTGM} we clearly have that $G_\eps(x,x)\to0$ as $x\to\partial\Omega$, so the restriction of $G_\eps$ to the set $\partial\Omega\times\overline\Omega\cup \overline\Omega\times \partial\Omega\cup\{(x,x)\ ;\ x\in\Omega\}$ is continuous. 

 Let $(x_0,y_0) \in \overline{\Omega}\times\overline{\Omega}$. We take $x \rightarrow x_0$ and $y \rightarrow y_0$ and we want to show that $G_\eps(x,y)\to G_\eps(x_0,y_0)$. We can assume without loss of generality that $x\neq y$ and $x,y\in\Omega$. We consider several cases depending on the location of $(x_0,y_0)$.

Assume first that $x_0\in \partial\Omega$ and $ y_0\in \overline{\Omega}$, with $x_0 \neq y_0$. By Theorem  \ref{lemmeGreenReguliere}, $$G_{\R^2}(x,y) + \gamma(x,y) = G(x,y) \tend{(x,y)}{(x_0,y_0)} G(x_0,y_0)= 0$$ since $x_0 \in \partial\Omega$. Moreover $$G_{\R^2}(x,y) \tend{(x,y)}{(x_0,y_0)}\frac{1}{2\pi} \ln |x_0-y_0|$$ so $$\gamma(x,y) \tend{(x,y)}{(x_0,y_0)}-\frac{1}{2\pi} \ln |x_0-y_0|.$$ We recall that $f_\eps$ is odd and continuous and thus $$G_\eps(x,y) = f_\eps(G_{\R^2}(x,y)) + f_\eps(\gamma(x,y)) \tend{(x,y)}{(x_0,y_0)} 0 = G_\eps(x_0,y_0).$$

Now assume $x_0 = y_0 \in \Omega$. Then 
$$G_{\R^2}(x,y) \tend{(x,y)}{(x_0,y_0)} - \infty$$
and
$$ \gamma(x,y)  \tend{(x,y)}{(x_0,y_0)} \tgm(x_0) $$
thus
$$ G_\eps(x,y)  \tend{(x,y)}{(x_0,y_0)} -L_\eps  + f_\eps(\tgm(x_0)) = G_\eps(x_0,y_0).$$
Finally, assume that $x_0 = y_0 \in \partial\Omega$. Then
$$G_{\R^2}(x,y) \tend{(x,y)}{(x_0,y_0)} - \infty$$
and by Corollary \ref{coroLimiteGamma} we have that
$$ \gamma(x,y)  \tend{(x,y)}{(x_0,y_0)} + \infty $$
thus
$$ G_\eps(x,y)  \tend{(x,y)}{(x_0,y_0)} -L_\eps +L_\eps  = 0 = G_\eps(x_0,y_0).$$

We conclude that $G_\eps$ is continuous. 

We prove now that $G_\eps$ is $C^1$ up to the boundary. Let us compute its gradient in the first variable for any $(x,y) \in \Omega\times\Omega, x \neq y$:
\begin{equation}\label{derGeps}
    \nabla_x G_\eps(x,y) = \nabla_x G_{\R^2}(x,y)f'_\eps(G_{\R^2}(x,y)) + \nabla_x \gamma(x,y)f'_\eps(\gamma(x,y)).
\end{equation}
Since $f_\eps'$ is compactly supported, $f'_\eps(G_{\R^2}(x,y)) = 0$ in a neighborhood of the diagonal  of $\overline{\Omega}\times\overline{\Omega}$. Similarly, if $x_0 \in \partial\Omega$ then $f'_\eps(\gamma(x,y)) = 0$ in a neighborhood of $(x_0,x_0)$ so $f'_\eps(\gamma(x,y))$ is smooth in a neighborhood of the diagonal of $\overline{\Omega}\times\overline{\Omega}$. 

Let now  $x_0 \in \partial\Omega$ and $y_0 \in \overline{\Omega}$ with $x_0 \neq y_0$. Consider $x\to x_0$, $y\to y_0$ where $x\neq y$ and $x,y\in\Omega$. By Theorem \ref{lemmeGreenReguliere}, $\nabla_x G(x,y)$ converges, so $\nabla_x \gamma(x,y)$ converges too, and thus all quantities involved in \eqref{derGeps} converge. We proved that $G_\eps \in C^1(\overline{\Omega}\times\overline{\Omega})$. The proof that $G_\eps \in C^2(\overline{\Omega}\times\overline{\Omega})$ follows along the same lines.

We now notice that $\nabla^\perp_x G_\eps(x,y)$ is tangent to the boundary when $(x,y) \in \partial\Omega\times\overline{\Omega}$, since $G_\eps(x,y)= 0$ when $x\in\partial\Omega$, for every $y\in\overline{\Omega}$.

We now prove the bounds stated in the lemma. Since $f_\eps$ is an odd Lipschitz map  with Lipschitz constant 1, we have that
\begin{equation*}
    \forall (a,b) \in \R^2, |f_\eps(x) + f_\eps(y)| = |f_\eps(x) - f_\eps(-y)| \le |x- (-y)| = |x+y| 
\end{equation*}
and therefore 
$$|G_\eps(x,y)|=|f_\eps(G_{\R^2}(x,y)) + f_\eps(\gamma(x,y))|
\leq |G_{\R^2}(x,y)+\gamma(x,y)| = |G(x,y)|.$$
Combining relation \eqref{greendecomp} with inequality \eqref{majDerGreen} we have that $ |\nabla_x \gamma(x,y) | \le \frac{C}{|x-y|}$ which gives that
\begin{equation*}
    |\nabla_x G_\eps(x,y)| \le \frac{1}{2\pi|x-y|} + \frac{C}{|x-y|}.
\end{equation*} 
Therefore $|\nabla_x G_\eps(x,y)| \le \frac{C}{|x-y|}$ and this completes the proof of the lemma.
\end{proof}

 Note that by construction, the implication \eqref{RegAvantTauEps} is true.
 We thus constructed a suitable regularized dynamics. 

\subsubsection*{Additional properties}
We need to establish that the estimates of Lemma \ref{lemmeMajScalGenerale} also hold true for the regularized dynamics.

\begin{lemme}\label{lemmeMajScalGeneralIndepEps}
We have that for any $\kappa < 1$
\begin{equation*}
     \iint_{\Omega\times\Omega}\frac{1}{|x-y|^\kappa}|\nabla_xG_\eps(x,y) \cdot \nabla^\bot\tgm_\eps(x)| \dd x\dd y < C.
\end{equation*}

\begin{equation*}
     \iint_{\Omega\times\Omega}\frac{1}{d(x,\partial\Omega)^\kappa}|\nabla_xG_\eps(x,y) \cdot \nabla^\bot\tgm_\eps(x)| \dd x\dd y < C.
\end{equation*}
where the constant $C$ doesn't depend on $\eps$.
\end{lemme}
\begin{proof}
One can check from \eqref{derGeps} that the following relation holds true
\begin{equation*}
        \nabla_x G_\eps(x,y) = \nabla_x G(x,y)f'_\eps(G_{\R^2}(x,y)) + \nabla_x \gamma(x,y)(f'_\eps(\gamma(x,y))-f'_\eps(G_{\R^2}(x,y))).
\end{equation*}
We use the expression of $\tgm_\eps$ and the previous relation to obtain that
\begin{multline*}
        \nabla_x G_\eps(x,y) \cdot \nabla^\bot\tgm_\eps(x) = \nabla_x G(x,y)\cdot \nabla_x^\perp\tgm(x) f'_\eps(\tgm(x))f'_\eps(G_{\R^2}(x,y)) \\+ \nabla_x \gamma(x,y)\cdot\nabla_x^\perp\tgm(x) f'_\eps(\tgm(x))(f'_\eps(\gamma(x,y))-f'_\eps(G_{\R^2}(x,y))) \\ \equiv A_{1,\eps}(x,y) + A_{2,\eps}(x,y).
\end{multline*}
Recalling that $|f_\eps'| \le 1$, we can apply directly Lemma \ref{lemmeMajScalGenerale} to the term $A_{1,\eps}$ to obtain that there exists a constant $C$ that doesn't depend on $\eps$ such that 
\begin{equation*}
     \iint_{\Omega\times\Omega}\frac{1}{d(x,\partial\Omega)^\kappa}|A_{1,\eps}(x,y)| \dd x\dd y < C
\end{equation*}
and
\begin{equation*}
     \iint_{\Omega\times\Omega}\frac{1}{|x-y|^\kappa}|A_{1,\eps}(x,y)| \dd x\dd y < C.
\end{equation*}
It remains to prove the same bounds for $A_{2,\eps}$. Let 
\begin{equation*}
    E =\{ (x,y) \in \Omega\times\Omega, \sp f'_\eps(\tgm(x))\big(f'_\eps(\gamma(x,y))-f'_\eps(G_{\R^2}(x,y))\big) \neq 0 \}.
\end{equation*}
Since $\big|f'_\eps(\tgm(x))\big(f'_\eps(\gamma(x,y))-f'_\eps(G_{\R^2}(x,y))\big)\big| \le 2$, we have that
\begin{equation}\label{majIntA2eps}
     \iint_{\Omega\times\Omega}p(x,y)|A_{2,\eps}(x,y)| \dd x\dd y \le \iint_{E}2p(x,y)|\nabla_x \gamma(x,y)||\nabla_x^\perp\tgm(x)| \dd x\dd y
\end{equation}
with $p(x,y) = \frac{1}{d(x,\partial\Omega)^\kappa}$ or $p(x,y) = \frac{1}{|x-y|^\kappa}$.

We now want to show that for every $(x,y) \in E$, we have that $d(x,\partial\Omega) \ge C \eps$. By construction of $f_\eps$, if $\tgm(x) > \frac{1}{2\pi}|\ln \eps|+ 1$ then $f_\eps'(\tgm(x)) = 0$. By construction of $E$ this means that for $\eps$ small enough such that $-(\frac{1}{2\pi}|\ln \eps|+ 1) < \min_{\Omega} \tgm$, for every $(x,y) \in E$, we have that $|\tgm(x)| \le \frac{1}{2\pi}|\ln \eps|+ 1$. Moreover relation \eqref{EquivTGM} gives that
\begin{equation*}
   -\ln d(x,\partial\Omega)-C_1 \le 2\pi|\tgm(x)|\leq |\ln\eps|+2\pi
\end{equation*} 
and therefore, provided $\eps < 1$,
\begin{equation*}
    d(x,\partial\Omega)  \ge \eps \exp\left(- 2\pi - C_1 \right) \equiv C_2 \eps.
\end{equation*}
Consequently $E \subset E_1$ with 
\begin{equation*}
    E_1 = \{(x,y)\in\Omega\times\Omega, d(x,\partial\Omega) \ge C_2\eps  \}.
\end{equation*}
Moreover we have that $f'_\eps(\gamma(x,y))-f'_\eps(G_{\R^2}(x,y)) = 0$ on the set 
\begin{equation*}
    E'=\big\{ (x,y)\in\Omega\times\Omega, |G_{\R^2}(x,y)| <  \frac{1}{2\pi}|\ln\eps| \text{ and }  |\gamma(x,y)| <  \frac{1}{2\pi}|\ln\eps|\big\}.
\end{equation*}
Since $E \subset (E')^c$, assuming that $\eps<\frac1{\diam\Omega}$ we have that $E \subset E_2 \cup E_3$ with
\begin{equation*}
    E_2 = \{ (x,y)\in\Omega\times\Omega, |x-y| \le \eps \}
\end{equation*} 
and
\begin{equation*}
    E_3 =  \big\{ (x,y)\in\Omega\times\Omega, |\gamma(x,y)| \ge  \frac{1}{2\pi}|\ln\eps| \big\}.
\end{equation*}

Using the fact that $ E \subset (E_1 \cap E_2) \cup (E_1\cap E_3)$ as well as relation \eqref{MajDerTGM} into relation \eqref{majIntA2eps} yields that
\begin{multline}\label{decompAB}
    \iint_{\Omega\times\Omega}p(x,y)|A_{2,\eps}(x,y)| \dd x\dd y \le \iint_{E_1 \cap E_2}p(x,y)B_\eps(x,y)\dd x \dd y \\ + \iint_{E_1 \cap E_3}p(x,y) B_\eps(x,y)\dd x \dd y
\end{multline}
with \begin{equation*}
    B_\eps(x,y) =  2|\nabla_x \gamma(x,y)| \frac{C}{ d(x,\partial\Omega)}.
\end{equation*}

We bound now the quantity $ \iint_{E_1 \cap E_2}p(x,y)B_\eps(x,y)\dd x \dd y$. Let $x\in\Omega$ be such that $d(x,\partial\Omega) \ge C_2 \eps$. We have that $\{ y\in\Omega, (x,y) \in E_2\} = D(x,\eps) \cap \Omega$. Since $\nabla_x\gamma$ is harmonic in both its variables, by the maximum principle we have that
\begin{equation*}
    \sup_{y \in D(x,\eps) \cap \Omega} |\nabla_x \gamma(x,y) | \le  \sup_{y \in  \partial(D(x,\eps) \cap \Omega)} |\nabla_x \gamma(x,y) |.
\end{equation*}
Since $\partial(D(x,\eps) \cap \Omega) \subset \{ y \in \C, |x-y| = \eps \} \cup \partial \Omega$, and $d(x,\partial\Omega) \ge C_2\eps$, we  have that for every $y\in \partial(D(x,\eps) \cap \Omega)$, there exists a constant $C$ independent of $x$ and $y$ such that $|x-y|> C\eps.$ Since $|\nabla_x \gamma(x,y)| \le \frac{C}{x-y}$, for every $y \in \partial(D(x,\eps) \cap \Omega)$ we have that $|\nabla_x \gamma(x,y)| \le \frac{C}{\eps}$. Therefore,
\begin{equation*}
    \sup_{y \in D(x,\eps) \cap \Omega} |\nabla_x \gamma(x,y) | \le  \frac{C}{\eps}
\end{equation*}
and thus in the case $p(x,y) = \frac{1}{d(x,\partial\Omega)^\kappa}$, 
\begin{align*}
   \iint_{E_1\cap E_2} p(x,y) B_\eps(x,y) \dd x\dd y & \le \frac{C}{\eps}\iint_{E_1\cap E_2}\frac{1}{d(x,\partial\Omega)^{1+\kappa}}\dd x \dd y \\ 
   & \le \frac{C}{\eps} 2\pi\eps ^2 \int_{ \{x, d(x,\partial\Omega) \ge C\eps\} } \frac{1}{d(x,\partial\Omega)^{1+\kappa}}\dd x \\ 
   & \le C \eps^{1-\kappa}\\
&\leq C
\end{align*}
where we used Lemma \ref{lemmeIneqIntermediaires}.

For the other expression of $p(x,y)$ we directly use that $|\nabla_x \gamma(x,y)| \le \frac{C}{|x-y|}$ as well as relation \eqref{MajDerTGM} to obtain that
\begin{align*}
   \iint_{E_1\cap E_2}p(x,y)B_\eps(x,y) \dd x\dd y 
&\le \iint_{E_1\cap E_2}\frac{C}{|x-y|^{1+\kappa}}\frac{1}{d(x,\partial\Omega)}\dd x \dd y \\
&\le \int _{d(x,\partial\Omega)\geq C_2\eps}
\frac{C}{d(x,\partial\Omega)}\Bigl(\int _{|x-y|\leq\eps}\frac{1}{|x-y|^{1+\kappa}} \dd y \Bigr)\dd x \\
&\le C \eps^{1-\kappa} \int _{d(x,\partial\Omega)\geq C_2\eps}
\frac{C}{d(x,\partial\Omega)}\dd x\\
&\le C \eps^{1-\kappa}|\ln(\eps)|\\
&\leq C.
\end{align*}
Therefore
\begin{equation*}
    \iint_{E_1\cap E_2} p(x,y) B_\eps(x,y) \dd x\dd y < C
\end{equation*}
for both choices of $p$.

Now we need to estimate the quantity $ \iint_{E_1 \cap E_3}p(x,y) B_\eps(x,y)\dd x \dd y$. We start by recalling that since $\gamma$ is smooth on $\Omega \times \Omega$ and symmetric, for $\eps$ small enough, the relation $(x,y) \in E_3$ implies that either $x \in V_j$ or $y\in V_j$ for an index $j$. By Lemma \ref{LemmeDiffGammabounded}, we know that on $V_j\times\Omega$ the map $\gamma - \gamma_{\Omega_j}$ is bounded by a constant $M$. We conclude that $\forall (x,y) \in E_3$, $\frac{1}{2\pi} |\ln \eps| - M \le |\gamma_{\Omega_j}(x,y)| $ and therefore using Lemma \ref{LemmeGammaMinoré} or Lemma \ref{LemmeGammaMinoréExterior} with $k=1$, we know that there exists a constant $C$ such that
\begin{equation*}
    |x-y| \le C \eps.
\end{equation*}
So $E_3$ is included in a domain of similar form to $E_2$ and we can reproduce the previous argument and conclude that there exists a constant $C$ independent of $\eps$ such that
\begin{equation*}
    \iint_{E_1\cap E_3} p(x,y) B_\eps(x,y) \dd x\dd y < C.
\end{equation*}
Recalling relation \eqref{decompAB}, we have proved that 
\begin{equation*}
     \iint_{\Omega\times\Omega}p(x,y)|A_{2,\eps}(x,y)| \dd x\dd y < C
\end{equation*}
which concludes the proof of the lemma.
\end{proof}

\begin{lemme}
There exists a constant $C$ independent of $\eps$ such that for every $1 \le j \le m$,
\begin{equation}\label{scalHarmoniqueEsp}
    |\nabla\tgm_\eps(x)\cdot \beta_j(x)| \le C.
\end{equation}
\end{lemme}
\begin{proof}
Recalling that $|f'_\eps|\le 1$, this lemma is a direct consequence of the fact that $\nabla\tgm_\eps(x) = \nabla\tgm(x)f_\eps'(\tgm(x))$ and of Corollary \ref{coroTgmOrtho}.
\end{proof}

\subsection{End of the proof of Theorem \ref{mainresult}}\label{sectionMainResult}
The end of the proof of Theorem \ref{mainresult} is largely inspired from the work previously done in \cite{marchioro1984vortex}.

Recall that $\Gamma = \{ X = (x_1,\ldots,x_N), d(X) > 0\}$ where 
\begin{equation*}
    d(X) = \min \left( \min_{i\neq j} |x_i - x_j|, \min_i \mathrm{d}(x_i,\partial\Omega) \right)\quad\forall X = (x_1,\ldots,x_N).
\end{equation*}

The aim of the Theorem \ref{mainresult} is to prove that $\tau(X) = \infty$ for $\lambda$-almost every $X$ in $\Omega^N$. Since $\lambda(\Omega^N \setminus \Gamma) = 0$, we can assume that $X\in \Gamma$. We denote by $S_tX = (x_1(t),\ldots,x_N(t))$ the maximal solution of the point vortex system \eqref{ptVortexDynamicCOmplete} with $(x_1(0),\ldots,x_N(0)) = X$, and by $S_t^\eps X = (x_1^\eps(t),\ldots,x_N^\eps(t))$ the global solution of equations \eqref{regularizedDynamic} which is the regularized dynamics with the same initial data $X$.

For any $X\in\Gamma$, we define $\tau_\eps(X)$ as the supremum of all times such that the system of relations
\begin{equation*}
    \begin{cases}
        |G_{\R^2}(x_i(t),x_j(t))| < \frac{1}{2\pi} |\ln \eps| 
        \\ |\tgm(x_i(t))| < \frac{1}{2\pi} |\ln \eps | 
        \\ |\gamma_{\Omega}(x_i(t),x_j(t))| < \frac{1}{2\pi} |\ln \eps|
    \end{cases}
\end{equation*}
are true for any $i\neq j$ and any $t \in [0, \tau_\eps(X)[$, with the convention that $\tau_\eps(X)=0$ if there exists no such $t$. By relations \eqref{RegAvantTauEps}, we know that for every $X\in\Gamma$, $S_t^\eps X = S_t X$ for every $t < \tau_\eps(X)$.

Let us introduce the function
\begin{equation}\label{DefF}
    F(r) = \exp(-\eta r),
\end{equation} with a constant $0 < \eta < 1$ that we will specify later. Let $\phi_\eps : \Gamma \rightarrow \R$ be defined by
\begin{equation}\label{defPHI}
    \phi_\eps(X) = \frac{1}{2} \sum_{i\neq j} F(G_\eps(x_i,x_j)) + \frac{1}{2}\sum_i F(-\tgm_\eps(x_i)),
\end{equation}
and $\Lambda_\eps : \Gamma \times \R_+ \rightarrow \R $ defined by
\begin{equation}\label{defLambda}
    \Lambda_\eps(X,t) = \der{}{t} \phi_\eps(S_t^\eps X).
\end{equation}
Since equations \eqref{regularizedDynamic} are autonomous, we have that
\begin{equation}\label{propLambda}
    \Lambda_\eps(X,t) = \Lambda_\eps(S_t^\eps X, 0).
\end{equation}
We now claim that the following proposition is true.
\begin{prop}\label{propTempsFiltrationCorrect} For every $X \in \Gamma$,
\begin{equation*}
\phi_\eps(S_{\tau_\eps(X)}^\eps X) \ge \frac{1}{2} F\left(\frac{1}{8\pi}\ln(\eps)\right) = \frac{1}{2} \eps^{-\frac{\eta}{8\pi}}.
\end{equation*}
\end{prop}
We delay the proof for the time being.
Let $\tau>0$ be a fixed time. By Proposition \ref{propTempsFiltrationCorrect}, 
$$\{ X\in\Gamma, \tau_\eps(X) \le \tau \} \subset \left\{ X\in\Gamma, \sup_{t\in [0,\tau]}\phi(S_t^\eps X) \ge \frac{1}{2} \eps^{-\frac{\eta}{8\pi}} \right\},$$
Therefore,
\begin{align*}
    \lambda( \{ X\in\Gamma, \tau(X)\le\tau \} ) & \le   \lambda( \{ X\in\Gamma, \tau_\eps(X) \le \tau \} ) \\
    & \le \lambda\left(\left\{ X\in\Gamma, \sup_{t\in [0,\tau]}\phi_\eps(S_t^\eps X) \ge \frac{1}{2} \eps^{-\frac{\eta}{8\pi}} \right\}\right) \\
    & \le 2\eps^{\frac{\eta}{8\pi}}\int_{\Gamma} \sup_{t\in [0,\tau]}\phi_\eps(S_t^\eps X) \dd \lambda(X).
\end{align*}
Recalling the definition \eqref{defLambda} of $\Lambda_\eps$ and relation \eqref{propLambda}, for every $t\in[0,\tau]$ we have that
\begin{equation*}
    \phi_\eps(S_t^\eps X) = \phi_\eps(X) + \int_0^{t} \Lambda_\eps(X,s)\dd s = \phi_\eps(X) + \int_0^t \Lambda_\eps(S_s^\eps X, 0)\dd s,\\
\end{equation*}
thus
\begin{equation*}
     \sup_{t\in [0,\tau]} \phi_\eps(S_t^\eps X) \le |\phi_\eps(X)| + \int_0^{\tau} |\Lambda_\eps(S_s^\eps X, 0)|\dd s.\\
\end{equation*}
Using Fubini-Tonelli's Theorem we have that
\begin{equation*}
     \int_\Gamma \sup_{t\in [0,\tau]}\phi_\eps(S_t^\eps X ) \dd \lambda(X) \le  \int_\Gamma |\phi_\eps(X)|\dd\lambda(X) + \int_0^\tau \int_{\Gamma}|\Lambda_\eps(S_s^\eps X,0)| \dd \lambda(X) \dd s.
\end{equation*}
Since the flow $S^\eps$ is Hamiltonian, it is area preserving (see  \cite[Corollary 1.10]{Arnold_Kozlov_Neishtadt}) and thus we have that for any $s \in \R_+$,
\begin{equation*}
    \int_\Omega  |\Lambda_\eps(S_s^\eps X,0) |\dd \lambda(X) = \int_{\Omega}| \Lambda_\eps(X,0)| \dd \lambda (X).
\end{equation*}
We will prove later that there exists a constant $A_0$ depending only on $\Omega$, $N$, $\eta$ and on the masses $(a_i)_i$, such that for every $\eps >0$ and $t\in \R$, we have that
\begin{equation}\label{majintphi}
    \int_{\Gamma}\phi_\eps(X)\dd \lambda(X) \le A_0,
\end{equation}
and
\begin{equation}\label{majderphi}
    \int_{\Gamma}\left|\Lambda_\eps(X,0)\right| \dd \lambda(X)\le A_0.
\end{equation} 
Relations \eqref{majintphi} and \eqref{majderphi} yield that
\begin{equation*}
     \int_\Gamma \sup_{t\in [0,\tau]}\phi_\eps(S_t^\eps(X)) \dd \lambda(X) \le  A_0 (1+\tau)
\end{equation*}
and therefore
\begin{equation*}
    \lambda( \{ X\in\Gamma, \tau(X) \le \tau \} ) \le 2\eps^{\frac{\eta}{8\pi}} A_0(1+\tau).
\end{equation*}
This being true for every $\eps >0$, and given that the left-hand side of the equation doesn't depend on $\eps$, letting $\eps\to0$ yields
\begin{equation*}
    \lambda( \{ X\in\Gamma, \tau(X) \le \tau \} ) = 0.
\end{equation*}
This is true for every time $\tau > 0$, and since 
\begin{equation*}
    \{ X\in\Gamma, \tau(X) < \infty \} = \ds \bigcup_{k \in \N^*}^\infty\{ X\in\Gamma, \tau(X) < k \},
\end{equation*} we have the desired result:
\begin{equation*}
    \lambda\{ X\in\Gamma, \tau(X) < \infty\} = 0.
\end{equation*}

\subsubsection*{Proof of Proposition \ref{propTempsFiltrationCorrect}}
We recall that for every $t \le \tau_\eps(X)$ and any $i\neq j$, we have that $G_\eps(x_i,x_j) = G(x_i,x_j)$, $\gamma_\eps(x_i,x_j) = \gamma(x_i,x_j)$ and $\tgm_\eps(x_i) = \tgm(x_i)$. Let $(x,y) \in \Omega$, $x\neq y$. We recall that at the time $t = \tau_\eps(X)$, there exist $i \neq j$ such that either
\begin{equation*}
    |G_{\R^2}(x_i^\eps(\tau_\eps(X)),x_j^\eps(\tau_\eps(X)))| = \frac{1}{2\pi} |\ln \eps| 
\end{equation*}
or
\begin{equation*}
    |\gamma_{\Omega}(x_i^\eps(\tau_\eps(X)),x_j^\eps(\tau_\eps(X)))| =\frac{1}{2\pi} |\ln \eps|
\end{equation*}
or
\begin{equation*}
    |\tgm(x_i^\eps(\tau_\eps(X)))| = \frac{1}{2\pi} |\ln \eps | .
\end{equation*}
Recalling the definition of $\Phi$ given by relation \eqref{defPHI}, and the fact that $F$ is positive, we have that
\begin{equation*}
\phi_\eps(S_{\tau_\eps(X)}^\eps X) \ge \max \left( \frac{1}{2}F(G_\eps(x_i^\eps(\tau_\eps(X)),x_j^\eps(\tau_\eps(X))), \frac{1}{2}F(-\tgm_\eps(x_i^\eps(\tau_\eps(X))))\right).
\end{equation*}
Therefore, since $F$ is decreasing, in order to prove Proposition \ref{propTempsFiltrationCorrect} it is enough to prove the following lemma.
\begin{lemme}
Let $(x,y)\in \Omega\times\Omega$. If one of the conditions
\begin{equation*}
    \begin{cases}
        |G_{\R^2}(x,y)| \ge \frac{1}{2\pi} |\ln \eps| 
        \\ |\gamma(x,y)| \ge\frac{1}{2\pi} |\ln \eps|
        \\ |\tgm(x)| \ge \frac{1}{2\pi} |\ln \eps | 
    \end{cases}
\end{equation*}
is true, then it implies that either
\begin{equation*}
        G(x,y) \le \frac{1}{8\pi}\ln \eps 
\end{equation*}
or
\begin{equation*}
        \tgm(x)\ge  -\frac{1}{8\pi} \ln  \eps.
\end{equation*}
\end{lemme}
\begin{proof}
Firstly, provided that $\eps$ is small enough such that $-\frac{1}{2\pi} |\ln \eps |< \min_\Omega \tgm$, which is possible by relation \eqref{tgmMinoree}, the relation $|\tgm(x)| \ge \frac{1}{2\pi} |\ln \eps | $ implies that $\tgm(x)\ge  -\frac{1}{2\pi} \ln  \eps \ge -\frac{1}{8\pi} \ln  \eps$. 

Secondly assume that $|\gamma(x,y)| \ge\frac{1}{2\pi} |\ln \eps|$. We recall the decomposition \eqref{decompOmegaVj}. Since $K$ is a compact set, the map $\gamma$ is bounded on $K\times K$. Therefore, provided $\eps$ is small enough such that $\frac{1}{2\pi} |\ln \eps |> \max_{K\times K} |\gamma|$, we have that the condition $ |\gamma(x,y)| \ge\frac{1}{2\pi} |\ln \eps|$ implies that there exists $0 \le j \le m$ such that either $x \in V_j$ or $y \in V_j$. By symmetry, we assume that $x \in V_j$. 

We recall Lemma \ref{LemmeDiffGammabounded} which states that the map $\gamma - \gamma_{\Omega_j}$ is bounded on $V_j \times \Omega$. Therefore, there exists a constant $M>0$ such that $|\gamma_{\Omega_j}(x,y)| \ge\frac{1}{2\pi} |\ln \eps| - M$. If $j=0$, we use Lemma \ref{LemmeGammaMinoré} with $k=1$ and $U=\Omega_j$, else we use Lemma \ref{LemmeGammaMinoréExterior} with $k=1$ and $\Pi = \Omega_j$, and $\mathcal{U} = \Omega$,  to obtain that $d(x,\partial\Omega) \le C \eps$. Therefore relation \eqref{EquivTGM} gives that
\begin{equation*}
    2\pi\tgm(x) \ge -\ln(C\eps) - C_1.
\end{equation*}
We deduce that there exists $\eps_0>0$ such that for every $\eps < \eps_0$ we have
\begin{equation*}
   \tgm(x) \ge -\frac{1}{8\pi}\ln(\eps).
\end{equation*}

Thirdly, assume that $|G_{\R^2}(x,y)| \ge \frac{1}{2\pi} |\ln \eps|$, which is equivalent to  $G_{\R^2}(x,y) \le \frac{1}{2\pi} \ln \eps$ provided that $\eps$ is smaller than $\frac{1}{\diam \Omega}$. Recalling relation \eqref{greendecomp}, then either $G(x,y) \le \frac{1}{4\pi}\ln \eps$ or $ -\gamma(x,y) \le\frac{1}{4\pi} \ln \eps$. The condition $G(x,y) \le \frac{1}{4\pi}\ln \eps$ naturally implies that $G(x,y) \le \frac{1}{8\pi}\ln \eps$. Assume now that $ \gamma(x,y) \ge-\frac{1}{4\pi} \ln \eps$. As in the second case, we use Lemma \ref{LemmeGammaMinoré} or Lemma \ref{LemmeGammaMinoréExterior} with $k=\frac{1}{2}$ to obtain that $d(x,\partial\Omega) \le C \sqrt{\eps}$, which leads by relation \eqref{EquivTGM} to 
\begin{equation*}
   \tgm(x) \ge -\frac{1}{8\pi}\ln(\eps),
\end{equation*}
for $\eps$ small enough.
The lemma is now proved, which concludes the proof of Proposition \ref{propTempsFiltrationCorrect}.
\end{proof}

\subsubsection*{Proof of relations \eqref{majintphi} and \eqref{majderphi}.}
We start by proving \eqref{majintphi}. Recalling the definitions of $\phi_\eps$ and $F$ given by relations \eqref{defPHI} and \eqref{DefF} we have that
\begin{equation*}    
\phi_\eps(X) = \frac{1}{2} \sum_{i\neq j} \exp(-\eta G_\eps(x_i,x_j)) + \frac{1}{2}\sum_i \exp(\eta\tgm_\eps(x_i)).
\end{equation*}
Since $\eta > 0$ and $\exp$ is an increasing function, relations \eqref{PropRegFunctions} yield that
\begin{equation*}    
|\phi_\eps(X)| \le  \sum_{i\neq j} \exp(\eta|G(x_i,x_j)|) +\sum_i \exp(\eta|\tgm(x_i)|).
\end{equation*}
Relation \eqref{EquivTGM} gives that $|\tgm(x_i)| \le -\frac{1}{2\pi}\ln d(x_i,\partial\Omega) + C_3$. Using also \eqref{majGreen}, we have that
\begin{equation*}
    |\phi_\eps(X)|  \le \frac{1}{2} \sum_{i\neq j} \exp(\eta C(1 + |\ln|x_i-x_j||)) + \frac{1}{2}\sum_i \exp(- \frac{\eta}{2\pi}\ln d(x_i,\partial\Omega) +\eta C_3 ).
\end{equation*}
Since $|\ln|x-y|| \le \max\{-\ln |x-y|,  \ln(\diam \Omega) \}$, we bound $|\ln|x-y|| \le -\ln |x-y| + C$ and thus
\begin{equation*}
     |\phi_\eps(X)| \le  \exp(\eta C_4)\left[ \frac{1}{2} \sum_{i\neq j} \left(\frac{1}{|x_i-x_j|^{\eta C}}\right)  + \frac{1}{2}\sum_i \frac{1}{d(x_i,\partial\Omega)^{\eta/2\pi}}\right].
\end{equation*}

Choosing $\eta < \min \{ \frac{2}{C}, 2\pi\}$, and noticing that the last expression doesn't depend on $\eps$, we obtain that $ \int_\Gamma \phi_\eps(X)\dd \lambda(X)$ is bounded independently of $\eps$. This proves relation \eqref{majintphi}.

We now want to prove relation \eqref{majderphi}. By the definition of $\Lambda_\eps$ given in relation \eqref{defLambda}, we have that
\begin{equation*}
    \int_{\Gamma}\left|\Lambda_\eps(X,0)\right| \dd \lambda(X) = \int_{\Gamma}\left|\der{}{t}[\phi_\eps(S_t^\eps X)]\bigl|_{t=0}\right| \dd \lambda(X).
\end{equation*} 
Therefore in order to prove relation \eqref{majderphi} we have to show that at time $t=0$, the quantity $  \der{}{t} \phi_\eps(S_t^\eps X) $ is bounded in $L^1(\Gamma)$ independently of $\eps$. Let us compute:
\begin{align*}
    \der{}{t} \phi_\eps(S_t^\eps X) & = \der{}{t} \left[ \frac{1}{2} \sum_{i\neq j} F(G_\eps(x_i^\eps(t),x_j^\eps(t))) + \frac{1}{2}\sum_i F(-\tgm_\eps(x_i^\eps(t)))\right] \\
    & = \frac{1}{2}\sum_{i\neq j}F'(G_\eps(x_i^\eps(t),x_j^\eps(t)))\nabla_xG_\eps(x_i^\eps(t),x_j^\eps(t)) \cdot \der{x_i^\eps(t)}{t}\\
    & \sp \sp + \frac{1}{2}\sum_{i\neq j}F'(G_\eps(x_i^\eps(t),x_j^\eps(t)))\nabla_yG_\eps(x_i^\eps(t),x_j^\eps(t)) \cdot \der{x_j^\eps(t)}{t}\\    
    & \sp \sp  - \frac{1}{2}\sum_i F'(-\tgm_\eps(x_i(t))) \nabla\tgm_\eps(x_i^\eps(t))\cdot  \der{x_i^\eps(t)}{t}.
\end{align*}
However, one can notice that since $G_\eps(x,y) = G_\eps(y,x)$ we have that $\nabla_y G_\eps(x,y) = \nabla_x G_\eps(y,x)$ and thus
\begin{align*}
    \der{}{t} \phi_\eps(S_t^\eps X) & = \sum_{i\neq j}F'(G_\eps(x_i^\eps(t),x_j^\eps(t)))\nabla_xG_\eps(x_i^\eps(t),x_j^\eps(t)) \cdot \der{x_i^\eps(t)}{t}\\    
    & \sp \sp  - \frac{1}{2}\sum_i F'(-\tgm_\eps(x_i^\eps(t))) \nabla\tgm_\eps(x_i^\eps(t))\cdot  \der{x_i^\eps(t)}{t}.
\end{align*}
We recall relation \eqref{regularizedDynamic}:
\begin{equation*}
    \der{x_i^\eps(t)}{t} = \sum_{\substack{k= 1\\ k\neq i}}^N \nabla_x^\bot G_\eps(x_i^\eps(t),x_k^\eps(t))a_k + \frac{1}{2} \nabla^\bot\tgm_\eps(x_i^\eps(t)) a_i + \sum_{k=1}^m c_k(t)\beta_k(x_i^\eps(t)).
\end{equation*}
Therefore,
\begin{align*}
    \der{}{t} \phi_\eps(S_t^\eps X) & = \sum_{\substack{i\neq j \\ k \neq i }}F'(G_\eps(x_i^\eps(t),x_j^\eps(t)))\nabla_xG_\eps(x_i^\eps(t),x_j^\eps(t)) \cdot \nabla_x^\bot G_\eps(x_i^\eps(t),x_k(t))a_k\\
    & \sp\sp + \sum_{i\neq j}F'(G_\eps(x_i^\eps(t),x_j^\eps(t)))\nabla_xG_\eps(x_i^\eps(t),x_j^\eps(t)) \cdot \frac{1}{2} \nabla^\bot\tgm_\eps(x_i^\eps(t)) a_i \\
    & \sp\sp + \sum_{i\neq j}F'(G_\eps(x_i^\eps(t),x_j^\eps(t)))\nabla_xG_\eps(x_i^\eps(t),x_j^\eps(t)) \cdot \sum_{k=1}^m c_k(t)\beta_k(x_i^\eps(t))\\
& \sp \sp  - \frac{1}{2}\sum_i F'(-\tgm_\eps(x_i^\eps(t))) \nabla\tgm_\eps(x_i^\eps(t))\cdot \sum_{\substack{k= 1\\ k\neq i}}^N \nabla_x^\bot G_\eps(x_i^\eps(t),x_k^\eps(t))a_k \\
& \sp \sp - \frac{1}{4}\sum_i F'(-\tgm_\eps(x_i^\eps(t))) \nabla\tgm_\eps(x_i^\eps(t))\cdot  \nabla^\bot\tgm_\eps(x_i^\eps(t)) a_i \\
& \sp \sp - \frac{1}{2}\sum_i F'(-\tgm_\eps(x_i^\eps(t))) \nabla\tgm_\eps(x_i^\eps(t))\cdot \sum_{k=1}^m c_k(t)\beta_k(x_i^\eps(t)) \\
& \equiv B_1(t) + B_2(t) + B_3(t) + B_4(t) + B_5(t) + B_6(t).
\end{align*}

We recall that $x_i^\eps(0) = x_i$, where $X = (x_1,\ldots,x_N)$. First of all, we observe that $B_5(t) = 0$. Notice that from  relations \eqref{PropRegFunctions} we have that
\begin{equation}\label{majDerGreenEps}
    |\nabla_xG_\eps(x,y)| \le \frac{C}{|x-y|}
\end{equation}
where the constant $C$ is independent of $\eps$.

The same estimates as at the beginning of the proof of relations \eqref{majintphi} and \eqref{majderphi} show that for any $X \in \Gamma$ we have that
\begin{equation}\label{majFprimeG}
    |F'(G_\eps(x_i,x_j))| \le \frac{C'}{|x_i-x_j|^{\eta C}}
\end{equation}
and
\begin{equation}\label{majFprimetgm}
    |F'(-\tgm_\eps(x_i))| \le \frac{C}{d(x_i,\partial\Omega)^{\eta/2\pi}}
\end{equation}
where we used that $0<\eta<1$.
Relation \eqref{majFprimetgm}, together with relation \eqref{scalHarmoniqueEsp}, yields
\begin{equation*}
    |B_6(0)| \le \sum_i\frac{C}{d(x_i,\partial\Omega)^{\eta/2\pi}}.
\end{equation*}
We also have that
\begin{equation*}
    |c_k(0)\beta_k(x_i)| \le  C,
\end{equation*}
and therefore using relations \eqref{majFprimeG} and \eqref{majDerGreenEps} we have that
\begin{equation*}
    |B_3(0)| \le  \sum_{i\neq j}\frac{C'}{|x_i-x_j|^{\eta C+1}}.
\end{equation*}
Both $B_3(0)$ and $B_6(0)$ are therefore bounded in $L^1(\Gamma)$ uniformly in $\eps$ provided that $\eta$ is small enough. 

Using relation \eqref{majFprimeG}, we have that
\begin{equation*}
    |B_2(0)| \le \sum_{i\neq j}\frac{C' }{|x_i-x_j|^{\eta C}}|\nabla_xG_\eps(x_i,x_j) \cdot  \nabla^\bot\tgm_\eps(x_i)| 
\end{equation*}
and thus Lemma \ref{lemmeMajScalGeneralIndepEps} implies that $B_2(0)$ is bounded in $L^1(\Gamma)$ uniformly in $\eps$ if $\eta$ is small enough. Similarly, using this time relation \eqref{majFprimetgm}, we have that
\begin{equation*}
    |B_4(0)| \le \sum_{i\neq k} \frac{C}{d(x_i,\partial\Omega)^{\eta/2\pi}} | \nabla\tgm_\eps(x_i)\cdot \nabla_x^\bot G_\eps(x_i,x_k)|
\end{equation*}
and once again, Lemma \ref{lemmeMajScalGeneralIndepEps} allows to conclude that $B_4(0)$ is bounded in $L^1(\Gamma)$ uniformly in $\eps$ if $\eta$ is small enough.

Finally, we bound $B_1(0)$ by noticing first that for $k=j$, the expression $\nabla_xG_\eps(x_i^\eps(t),x_j^\eps(t)) \cdot \nabla^\bot_x G_\eps(x_i^\eps(t),x_k(t))$ vanishes, and therefore using relation \eqref{majFprimeG} and the relation \eqref{majDerGreenEps} we obtain that
\begin{equation*}
    |B_1(0)| \le\sum_{\substack{i\neq j \\ k \neq i \\ k \neq j}}\frac{C' }{|x_i-x_k||x_i-x_j|^{\eta C+1}}.
\end{equation*}
This term is bounded in $L^1(\Gamma)$ uniformly in $\eps$ when $\eta$ is small enough. This concludes the proof of relation \eqref{majderphi}. The proof of Theorem \ref{mainresult} is completed.

\paragraph{Acknowledgments.}
The author wishes to acknowledge very helpful discussions with Paolo Buttà and Carlo Marchioro on the subject of this paper. He also wishes to thank the Dipartimento di Matematica, Sapienza Università di Roma for its hospitality. Finally, he wishes to thank Dragoș Iftimie for his precious help and numerous advices.

\bibliographystyle{plain}
\bibliography{base}

\bigskip

\noindent \textbf{Martin Donati}: Université de Lyon, CNRS, Université Lyon 1, Institut Camille Jordan, 43 bd. du 11 novembre,
Villeurbanne Cedex F-69622, France.
\\
Email: donati@math.univ-lyon1.fr

\end{document}